\documentclass{article}

\usepackage{epsfig,graphicx,float,subcaption,diagbox}
\usepackage{amsmath}
\usepackage{authblk}

\newcommand{\be}{\begin{equation}}
\newcommand{\ee}{\end{equation}}

\title{Spectral solutions of PDEs on networks }

\author[1]{\normalsize{M. Brio}\thanks{brio@math.arizona.edu}}
\author[2]{\normalsize{J. G. Caputo }\thanks{caputo@insa-rouen.fr}}
\author[1]{\normalsize{H. Kravitz} \thanks{hkravitz@math.arizona.edu}}

\affil[1]{Department of Mathematics, University of Arizona, Tucson, Arizona 85721, USA.}
\affil[2]{Laboratoire de Math\'ematiques, INSA de Rouen Normandie\\ 76801 Saint-Etienne du Rouvray, France.}

\date{\ }

\begin{document}
\maketitle
\vspace{-1cm}

\date{\ }

\begin{abstract}
{
To solve linear PDEs on metric graphs with standard
coupling conditions (continuity and Kirchhoff's law), we 
develop and compare a spectral,
a second-order finite difference, and a discontinuous Galerkin
method. 
The spectral method yields eigenvalues and eigenvectors 
of arbitary order with machine precision and 
converges exponentially. These eigenvectors provide a 
Fourier-like basis on which to expand the solution;  
however, more complex coupling conditions require additional research.
The discontinuous Galerkin method provides approximations of arbitrary 
polynomial order; however computing high-order eigenvalues accurately 
requires the 
respective eigenvector to be well-resolved. The method allows 
arbitrary non-Kirchhoff flux conditions and requires 
special penalty terms at the vertices to enforce continuity of 
the solutions. 
For the finite difference method, the standard one-sided second-order 
finite difference stencil reduces the accuracy of the vertex solution to $
O(h^{3/2})$.
To preserve overall second-order accuracy,
we used ghost cells for each edge. 
For all three methods we provide 
the implementation details, their validation, and examples illustrating 
their performance for the eigenproblem, Poisson equation, and the 
wave equation. } \\
 Keywords:
{partial differential equations, metric graphs, spectral method, finite
difference, discontinuous Galerkin}
\end{abstract}

\date{\ }

\tableofcontents

\section{Introduction }

Partial differential equations (PDEs) on networks arise in many physical applications 
such as gas and water networks \cite{herty10}  \cite{water},
electromechanical waves in a transmission grid \cite{kundur},
air traffic control \cite{airTrafic}, and random 
nanofibre lasers \cite{Gaio}, to name a few. The underlying 
mathematical model consists of a metric graph: a finite set of vertices 
connected by arcs (oriented edges) on which a metric is assigned. At 
the vertices we can have different coupling conditions. The simplest 
assumes continuity of the field and zero total gradient at the 
vertices (Kirchhoff's law). The standard one-dimensional Laplacian 
together with 
these boundary conditions results in a generalised Laplacian
and associated Helmoltz eigenvalue problem. It can be shown that with 
these coupling conditions (continuity and Kirchhoff's law) the problem 
is self-adjoint, see \cite{berkolaiko17}, yielding real
eigenvalues and orthogonal eigenvectors. The eigenvectors form 
a complete basis of the appropriate set
of square integrable functions on the graph.

This spectral framework plays a key role for linear PDEs as we review and apply it in the present article.
Using it, the wave equation on a metric graph is treated exactly
as the one-dimensional wave equation on a one-dimensional finite interval with simple 
boundary conditions, see for example \cite{Hildebrand}.
Due to these strong geometrical properties, metric graphs have been studied 
extensively for the  Schr\"odinger
operator (or Helmholtz), see \cite{Gnutzmann} for a review
and the recent book \cite{sym19}
for research on nonlinear PDEs where the self-adjoint framework 
is not applied. The graph aspect (effect of the large network 
structure on the PDE solutions) is often overlooked in these studies 
and remains an undeveloped open area of research.
There are few studies of computational methods for PDEs on metric 
graphs. Exceptions are studies for water and gas networks in 
the engineering context, see \cite{water} for water networks 
and \cite{herty10} for gas networks. Note also our article
\cite{cd18} on the sine-Gordon equation. Most 
methods used there rely 
on finite difference (FD) spatial discretisations, while a recent 
study \cite{Arioli} introduced a
finite element method (FEM)  that is first-order at the vertices with 
second-order approximation inside the edge. The method was applied 
to the computations of the quantum graph spectrum, the solution 
of the elliptic equation, and the time evolution of the diffusion equation.

In this article, we present a systematic procedure to compute 
eigenvalues and eigenvectors of arbitrary order
for a general metric graph and use this spectral framework
to compute solutions of linear PDEs on graphs. We compare this spectral 
method to a second-order FD discretisation using centered finite differences at 
the inner points of the edges and also second-order approximation 
to the solution at the vertices. To generalise the finite
element approach of \cite{Arioli}, we introduce a Discontinuous 
Galerkin (DG) method of arbitrary polynomial order that allows 
for exact enforcement of Kirchhoff's law at the vertices while 
the continuity of the solution is enforced via penalty 
term(s) added to the weak formulation. These three methods are 
compared for the following PDE problems: a generalised Helmholtz, 
the Poisson equation and a damped wave equation (telegrapher’s equation).

We show that the spectral method is superior for the Helmholtz 
problem since eigenvalues/eigenvectors of arbitrary order can be 
computed without resolving the spatial structure of the 
respective eigenvectors. The FD and the DG methods lose accuracy as 
the order of the eigenvalue is increased because the spatial structure 
of the corresponding eigenvector is not resolved. We derive a second-order FD method for the solution at the vertices using ghost points 
for each edge sharing the same vertex to enforce the Kirchhoff 
conditions. The standard second-order FD approximations to the 
derivatives in the Kirchhoff equation result in a 
reduction of the accuracy of the solutions at the vertices to $O(h^{3/2})$.

The DG of high polynomial order gives very accurate eigenvalues 
for the resolved modes without any search procedures needed for 
nonlinear eigenvalue solvers of the spectral method. On the other hand, 
the solution of the Poisson equation for the DG with common penalty 
terms to enforce the continuity of the solution is inaccurate due to 
the large underlying discontinuous solution space unless unresolved 
modes are filtered out.

For linear evolution PDEs, the spectral components of the solution 
are fixed by the initial conditions, and the corresponding mode amplitudes 
decrease exponentially, similar to the one-dimensional linear evolution PDEs on 
a finite interval, see \cite{Hildebrand}.
Using a one-dimensional Fourier Transform 
on the initial data of the linear evolution PDEs and on the right hand 
side of the Poisson equation, an estimate of the modes that are necessary 
to be resolved for the required accuracy is obtained.

The article is organised as follows. The statement of the problem and
a brief review of the background information on metric quantum graphs are 
provided in Section 2. In Section 3, we describe the numerical spectral 
method. Section 4 presents the FD and DG methods, emphasising the implementation details. These three 
methods are applied to the Helmholtz, Poisson, and the telegrapher’s 
equations, and the results are discussed in Section 5. We conclude 
with the discussion of the obtained results and ideas for 
future work in Section 6.

\section{Wave equation on the graph}

We consider a finite metric graph with $n$ vertices
connected by $m$ edges of length $l_j, \; j=1:m.$
Each edge is parameterised by its length and is oriented arbitrarily
from the origin vertex $x=0$ to the end vertex $x=l_j$. Following
standard graph theory, we call these oriented edges arcs. We recall the definition of the {\it degree} of a vertex:  it is 
the number of edges connected to it.

On this graph, we define the vector component wave equation
\be\label{vwave} U_{tt} -{\tilde \Delta } U=0,\ee
where $U \equiv (u_1,u_2, \dots,u_m)^T$. Each component
satisfies the one-dimensional wave equation inside the respective arc,  
\be\label{wave}
{u_j}_{tt} -{u_j}_{xx} =0, ~~~  j=1, 2, \dots ,m \ee
In addition, at the vertices the solution should be continuous and also satisfy the 
Kirchhoff flux conditions at each vertex of degree $d$

\be\label{kirchof}
\sum_{j=1}^{d} \; {u_j}_{x} =0,
 \ee
where $\displaystyle {u_j}_{x}$ represent the outgoing fluxes for arc $j$ emanating from the vertex $p$.

Consider equation (\ref{vwave}). Since the problem is linear, we can separate time and space
and assume a harmonic solution $U(x,t)= e^{ikt} \: V(x)$.
We then get a Helmholtz or Schr\"{o}dinger eigenproblem for $V$ on the graph
\be\label{helm}
-{\tilde \Delta } V = k^2 V ,\ee
together with the above-mentioned coupling conditions at the vertices
and where ${\tilde \Delta }$ is the generalised Laplacian, i.e. the
standard Laplacian on the arcs together with the coupling conditions at
the vertices.
This generalised eigenvalue problem admits an inner product
obtained from the standard inner product on $L_2$ space, see \cite{solomon}. We have
\be\label{inproduct}
< f |g > \equiv \sum_{arc ~~j} < f_j |g_j >,~~~~< f_j |g_j >=\int_{0}^{l_j} f_j(x) g_j(x) dx . \ee
A  solution in terms of Fourier harmonics on each branch $j$ of length $l_j$ is 
\be\label{simple}
v_j(x) = A_j \sin k x + B_j \cos k x . \ee
It has been shown that, for the standard coupling conditions used here
(continuity and Kirchhoff), the eigenvectors $V^i$ form a complete 
orthogonal basis of the Cartesian product 
$L_2([0,l_1]) \times L_2([0,l_2])  \dots \times L_2([0,l_m]) $, see \cite{solomon}.
Writing down the coupling conditions at each vertex, one obtains a 
homogeneous linear system whose $k$-dependent matrix is 
singular at the eigenvalues.

Using solution (\ref{simple}) on each arc with unknown coefficients $A_j$
and $B_j$, the coupling conditions at each vertex yield 
the homogeneous system
\be\label{Meqn}
M(k) X =0, \ee
of $2m$ equations for the vector of $2m$ unknown arc amplitudes  
$$\displaystyle X=(A_1,B_1, A_2, B_2, \dots A_m, B_m)^T  . $$ 
The matrix $M(k)$ is singular at the eigenvalues $-k^2$. We call
these $k$-values resonant frequencies.
A practical, robust computational algorithm for the computation of these
eigenvalues and eigenvectors is presented in the next section.

For each resonant frequency $k_q$, the eigenvectors $V^q$ are determined from the null space of matrix $M(k_q)$. They can then be written as
\be\label{Vi}
V^q = 
\begin{pmatrix} A^q_1 \sin k_q x + B^q_1 \cos k_q x \cr 
                A^q_2 \sin k_q x + B^q_2 \cos k_q x \cr
                \dots  \cr
                A^q_m \sin k_q x + B^q_m \cos k_q x \cr
\end{pmatrix}
\ee
They can be normalised using the scalar product defined above. We have
\be\label{vivi}\lVert V^q \rVert^2 =< V^q V^q> = \sum_{j=1}^m < V^q_j V^q_j>, 
\ee
where $V^q_j = A^q_j \sin k_q x + B^q_j \cos k_q x$ and $< V^q_j V^q_j>$
is the standard scalar product on $L_2([0,l_j])$ . This defines
a broken L$_2$ norm or graph norm. 
The scalar product $< V^q_j V^q_j>$ can be computed explicitly 
\be\label{vpvp}
< V^q_j V^q_j>= \left ( {A^q_j}^2 + {B^q_j}^2 \right ) {l_j \over 2}
+ {\sin {2 k_q l_j} \over 4 k_q} \left ( -{A^q_j}^2 + {B^q_j}^2 \right )
+ A^q_j B^q_j { 1 - \cos {2 k_q l_j} \over 2 k_q}.\ee

Once the eigenvalue problem is solved, 
one can proceed with the spectral solution of 
the time-dependent problem, exactly as for the one-dimensional wave equation
on an interval. 
For that, expand the solution of the wave equation on the
graph (\ref{vwave}) using the eigenvectors, 
\be\label{spec1}
U(x,t) = \sum_{q=1}^\infty a_q(t) V^q ,\ee
and obtain a simplified description of the dynamics in terms of the amplitudes
$a_q$.

Let us consider the initial value problem for equation (\ref{vwave}) with
$U(t=0)=U_0,~~U_t(t=0)=d U_0$. Plugging the expansion (\ref{spec1}) 
into (\ref{vwave}) and projecting on the eigenvector $V^q$ we get the
following amplitude equation
\be\label{aitt}
{d^2 {a_q} \over d t^2}  + k_q^2 {a_q}=0. \ee
The initial conditions for $a_q$ are obtained as usual by projecting $U(t=0)$
and $U_t(t=0)$ on $V^q$.
We have
\be\label{iniai}
a_q(t=0) = < U_0 V^q >,~~  {d a_q \over dt} (t=0) = < dU_0~ V^q > . \ee

For example, assuming an initial condition has only support on the first arc, 
\be\label{lpsi0}
U_0 =
\begin{pmatrix}  g(x) \cr 
                   0 \cr
                \vdots  \cr
                0 \cr 
\end{pmatrix}, ~~
dU_0 =
\begin{pmatrix}  g'(x) \cr 
                   0 \cr
                \vdots  \cr
                0 \cr 
\end{pmatrix}, ~~
\ee
we get
\be\label{ailpsi0} 
a_q(t=0) = < g(x) V^q_1(x) >, ~~{d a_q \over dt} (t=0) = <g'(x) V^q_1(x) > .\ee

\section{Spectral algorithm }

The procedure to compute the eigenvalues and eigenvectors for arbitrary graphs
involves the following steps. 

\begin{enumerate}
\item Form matrix $M(k)$ from the continuity and Kirchhoff linear equations
using symbolic manipulations described below.
\item Plot inverse condition number $\text{rcond(M(k))}$ as function of $k$ to estimate 
graphically the lower bound of the spacing $\Delta k$ between the consecutive values of $k$.  
\item Bracket each resonant frequency by splitting
the range into smaller sub-intervals.
\item Apply a line minimisation algorithm to the 
function $\text{rcond}(M(k))$ on each sub-interval to estimate each $k$
within a chosen tolerance. 
\item For each resonant frequency, perform a singular value decomposition (SVD) 
of the matrix $M(k)$ to determine its null space. This will provide the 
corresponding eigenvector $V$ in (\ref{Vi}).
\end{enumerate}

\subsection{Generation of the matrix $M(k)$}

To generate the matrix $M(k)$, we use the incidence matrix with rows and columns 
corresponding to the vertices and arcs of the graph. Each column then has two 
nonzero entries, $1$ and $-1$ reflecting the orientation chosen for this arc. 

On each arc, $e_i$, for each given resonant frequency $k$, the solution 
has the form $u_i(x)=A_i\sin(k x) + B_i\cos(k x)$. 
The continuity condition between arcs $e_i$ and $e_j$ yields
\be
 A_i\sin(k x_*) + B_i\cos(k x_*) - \left[ A_j\sin(k x_*) + B_j\cos(k x_*)\right] =0
 \ee
where
 $$x_*=
\begin{cases}
0 \text{, if arc } e_r \text{ is outgoing}\\
l_r\text{, if arc } e_r \text{ is incoming}
\end{cases}
$$

To implement the continuity junction condition, our algorithm first finds all unique arc pair combinations at each vertex and saves them to a cell array $C$ (one cell for each vertex containing a list of pairs). 
For simplicity, the arc with the lowest arc number is fixed and the continuity conditions are written for this arc with every other arc. For a vertex of degree $d$, there will be $d-1$ continuity equations.

For each vertex, the algorithm loops through all of the continuity equations and places respective sines and cosines for each pair $ A_i, B_i$ into a row of matrix $M(k)$. In particular,  each continuity equation will place the sine term for the $i$th arc in the $(2i-1)$th column, and the cosine term in the $2i$th column schematically shown below, 

$$\begin{bmatrix}

\sin(k x_*) & \cos(k x_*) & -\sin(k x_*)& -\cos (k x_*) & \cdots & \vdots \\
\cos(k x_*) & -\sin(k x_*) & \cos(k x_*)& -\sin (k x_*) & \cdots & \vdots \\
\vdots & \vdots & \vdots& \vdots & \vdots & \vdots \\
\cdots & \cdots & \cdots& \cdots & \cdots & \cdots 
\end{bmatrix}
\begin{bmatrix}
A_1 \\
B_1 \\
A_2 \\
B_2 \\ 
\vdots \\
\end{bmatrix}
=0
$$

The Kirchhoff flux vertex condition is implemented at each vertex by summing the derivatives. Note that the $k$ term in front of the sine and cosine will cancel when the sum of the derivatives are set equal to zero and can be omitted. The equation for the flux condition is

\be \sum_{i} A_i \cos(k x_*) - B_i \sin(k x_*)=0.  \ee 

This equation takes up one row of the matrix for each vertex, thus resulting in $d$ equations for a vertex of degree $d$, with $d-1$ equations from the continuity constraint and one from the flux condition. 
The result of applying these conditions is the system of equations (\ref{Meqn}) 


Note that languages with symbolic capabilities such as Matlab, 
Mathematica, etc. have a command that converts linear equations 
into matrix form using the list of equations and the list of unknowns 
($\displaystyle X=(A_1,B_1, A_2, B_2, \dots A_m, B_m)^T$) as inputs.  

\subsection{Finding the resonant $k$s} 
 
\begin{figure}[H]
\centerline{
\epsfig{file=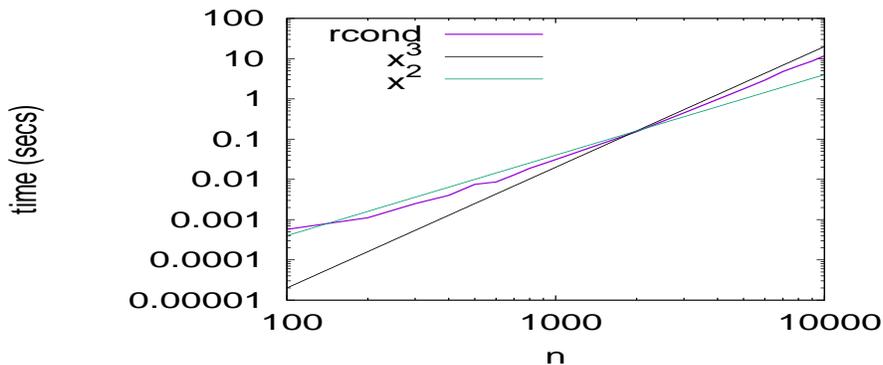,height=5cm,width=12 cm,angle=0}
}
\caption{Plot of time in secs as a function of $m$ for rcond.}
\label{rcond}
\end{figure}

Once the matrix $M(k)$ is formed, one estimates its inverse condition
number $r(k)$, 
\be\label{rk} r(k) = 
{1 \over \parallel M \parallel_1 \parallel M^{-1} \parallel_1} , \ee
using the rcond Matlab function, whose
complexity is 
$  O( 4 m^2)$, see \cite{moler}.
In fact, the computation of rcond involves an LU factorisation and its
complexity can be estimated using Matlab. We find it to be 
$O(m^2)$ for $m < 1000$ and $O(m^3)$ for larger $m$, see Fig. \ref{rcond}.

Note that the resonant frequencies are calculated via the robust 
inverse condition number estimator $\text{rcond(M(k))}$
instead of a determinant computation that 
requires arbitrary precision arithmetic for an accurate evaluation. 
At the resonant $k$s , $M(k)$ is singular and $M^{-1}$
ceases to exist. In practice $\parallel M^{-1} \parallel \to \infty$
so that $r \to 0$.

To find $k^*$ such that $r(k^*)=0$, we implemented a gradient-less
line minimisation algorithm, see \cite{linesearch}; this is similar
to a bisection. This gave $k^*$ up to machine precision.

The eigenvector $V^*$ associated to the eigenvalue $-{k^*}^2$
of the generalised Laplacian (\ref{helm}) on the graph can be found using the
SVD of $M(k^*)$ 
$$M(k^*) = U \Sigma V^T  . $$ 
Since $M(k^*)$ is singular, it has at least one zero singular value. The column vector of $V$ associated to the zero singular value spans the
kernel of $M(k^*)$.

\section{Finite difference and discontinuous Galerkin algorithms}

We present here finite difference discretisations and introduce
a discontinuous Galerkin approximation.

\subsection{Finite difference differentiation matrix on the graph} 

Here we describe a second-order finite difference method for both 
the inner arc discretisation and the vertex conditions. The approach can 
be generalised to higher-order finite differences.

Each arc $e_i$ is discretised into $N_i+1$ equally-spaced points enumerated from $j=0$ to $j=N_i$, with the number of points varying for each edge. $j=0$ and $j=N_i$ represent the starting and ending vertex of the edge respectively.
Each arc has its own uniform spatial mesh spacing $\Delta x_i$. We use a centered, explicit second-order scheme for the second derivative in space,  \\
\be
u_{xx} = \frac{u_{j+1,i} - 2u_{j,i} + u_{j-1,i}}{\Delta x_i^2} + {O} \left( \Delta x_i^2 \right)
\ee
where $u_{j,i}  \approx u(x_j)$, the solution at $x_j$ on a particular edge $e_i$ for inner points $j=1, 2, \dots, N_i-1$.

For the Poisson equation we solve the following FD equation at the inner points
\be \frac{u_{j+1,i} - 2u_{j,i} + u_{j-1,i}}{\Delta x_i^2} = f(x_j). \ee

To implement the vertex conditions, we label arcs adjacent to the vertex as $c=1, 2, \dots, d$. To enforce the continuity condition we place a node exactly at the vertex and denote the solution value as $u_0$. The center vertex $u_0$ 
is shared by each adjacent arc via the continuity condition $u_{c,0}=u_0 $. To implement the Kirchhoff flux condition, we use a centered, second-order scheme for the first derivative by extending each arc and adding a ghost point next to the vertex at $x=x_{-1,c}$ for each adjacent arc $e_c$. The derivative is taken in the outgoing direction from the center vertex.
For $j=0$ we have 
\be 0 = \sum_{c=1}^d u_x(x_0) = \sum_{c=1}^d \frac{u_{1,c} - u_{-1,c}}{2 \Delta x_c}  + {O}(\Delta x_1^2 + \Delta x_2^2 + \cdots + \Delta x_d^2) \ee
Here we say $j=1$ is the point adjacent to the vertex regardless of edge orientation.

We then apply the FD scheme for the PDE. For the Poisson equation we have
\be \frac{u_{1,c} - 2 u_0 + u_{-1,c}}{\Delta x_c^2} = f(x_0) \ee
This equation can be solved for $u_{-1,c}$ and substituted into the Kirchhoff flux equation to eliminate the ghost point. 
\be 0 =  \sum_{c=1}^d \frac{u_{1,c} - \left[ \Delta x_c^2 f(x_0) +2 u_0 - u_{1,c} \right]}{2 \Delta x_c}   \ee

This results in the following equation for 
$u_0$ and its neighboring values for placement in the finite difference matrix.
\be 2 \frac{ \sum_{c=1}^d \frac{u{1,c}}{\Delta x_c}  - u_0  \sum_{c=1}^d \frac{1}{\Delta x_c}}{ \sum_{c=1}^d \Delta x_c} =f(x_0) \ee

\begin{figure}[H]
\centerline{
\epsfig{file=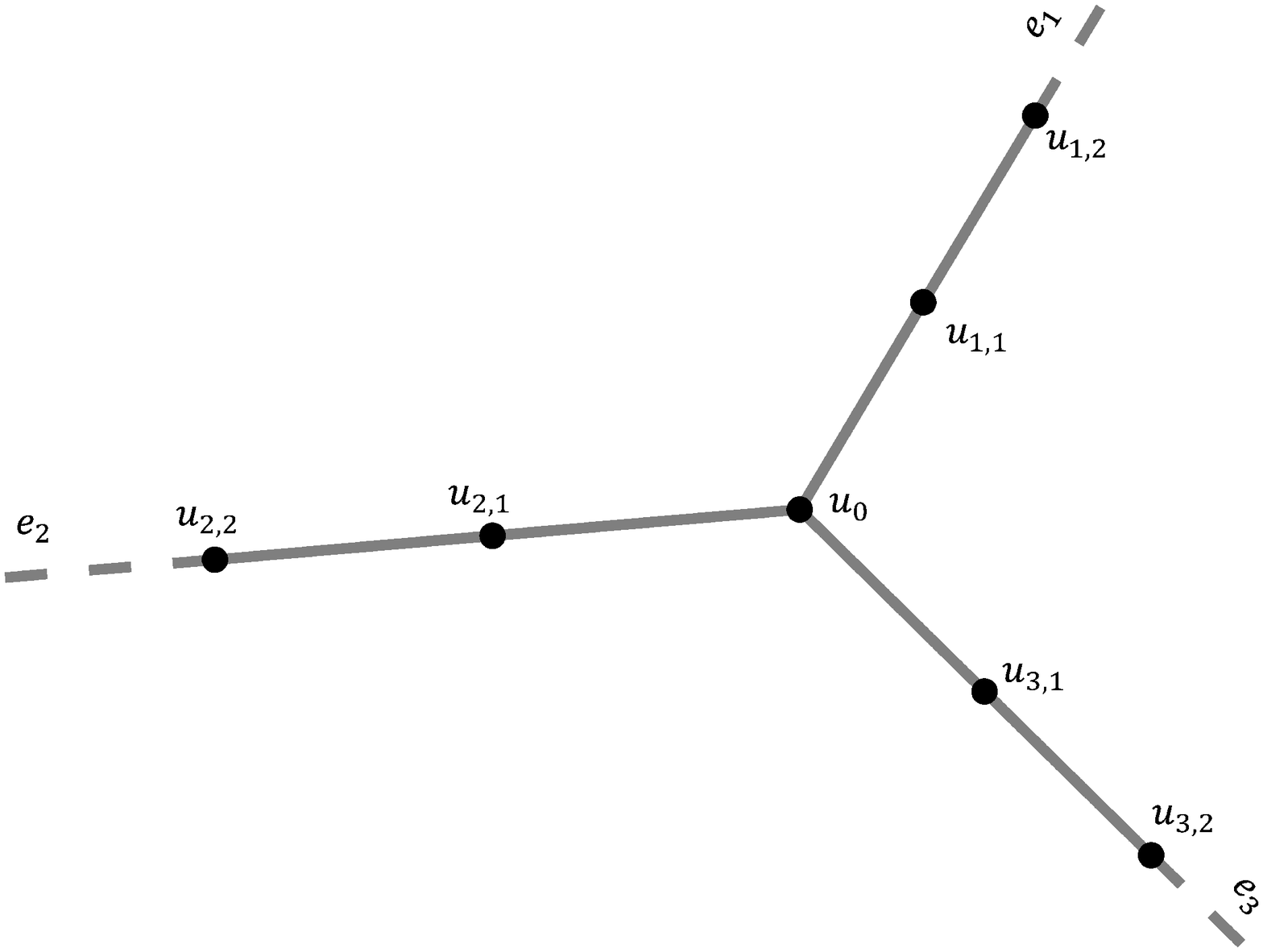,height=5cm,width=10 cm,angle=0}}
\caption{Finite difference discretisation for a vertex of degree 3.}
\label{u0}
\end{figure}

It is important that the orders of the approximation method at 
the vertices and in the bulk of the arcs be the same.  \cite{Arioli} use a piecewise linear basis and enforce continuity of the solution by 
placing a node at the vertex. 
This renders the approximation first-order, 
despite the inner approximation being second order. In our FD
scheme, the orders at the vertex and in the bulk of the arcs are
the same.

\subsection{Discontinuous Galerkin differentiation matrix on the graph}

The Discontinuous Galerkin Finite Element Method (DG) is a multi-domain 
method where the
solution is approximated by polynomials on each subdomain with 
appropriate interface conditions. 
This allows flexibility in choosing the degree of the polynomials 
used on the different subdomains and explicitly introduces the 
interface fluxes in the weak formulation of the problem containing 
second-order spatial partial derivatives.

In this article, we use 
orthogonal Legendre polynomials. The Kirchhoff flux  
conditions are applied explicitly in the weak formulation 
(natural interface boundary conditions). The continuity of the solution 
is enforced via a penalty term. In contrast, in the standard finite element 
approximation introduced in \cite{Arioli},
the continuity is enforced by placing a node at each vertex. This 
resulted in a first-order approximation of the Kirchhoff conditions 
even though a second-order (piecewise linear) approximation was used 
inside each arc.

We verify that the error estimate for the Discontinuous Galerkin
in the $L_2$-norm is $O(h^{p+1})$ where $h$ is the size of 
the subdomain and $p$ is the degree of the polynomial used on 
the different subdomains. Note that DG is naturally suitable for 
an $h-p$ adaptive strategy to determine the most efficient 
combination between $h$ and $p$ for the problem at hand.

We illustrate the implementation of the DG for the
quantum graph eigenvalue problem. 
Each arc is split into 
its own number of intervals, $[x_j,x_{j+1}]$,  and 
equation 
\be  u_{xx}=-k^2 u \ee 
is rewritten in a weak form for numerical approximation $u(x)$ using 
test function $\phi(x)$
\be 
 \int\limits_{x_j}^{x_{j+1}} u_{xx} \phi \; {\rm d}x =-k^2 \int\limits_{x_j}^{x_{j+1}} u_x  \phi \; {\rm d}x . 
\ee
The second derivative term is rewritten in equivalent ultra-weak form after 
applying three integration by parts,  see \cite{Chen}, 
\be \label{dg1}
 \int\limits_{x_j}^{x_{j+1}} u_{xx}\phi \; {\rm d}x =-\int\limits_{x_j}^{x_{j+1}} u_x \phi_x \; {\rm d}x + u_x \phi  \Bigr|_{x_{j}}^{x_{j+1}} + (u-u) \phi_x  \Bigr|_{x_{j}}^{x_{j+1}}.
\ee 

Among numerous penalty methods for the DG method applied to 
the second-order spacial derivatives 
\cite{Chen}, 
we chose the following approach that can be generalised to the graph vertices to enforce continuity there, and also allows for applying both Kirchhoff and non-Kirchhoff scattering flux conditions at the vertices to be explored in the future. 

For the inner intervals (endpoints do not include the vertices) the communication between neighboring intervals is introduced through both boundary terms that are evaluated as follows, 
\be 
{\widetilde u}_x \phi  \Bigr|_{x_{j}}^{x_{j+1}} + (u-{\widehat u}) \phi_x  \Bigr|_{x_{j}}^{x_{j+1}} ,
\ee 
where the $\widehat u$ denotes the arithmetic average across the jumps, while ${\widetilde u}_x $ has an additional penalty term to enforce continuity of the solution across the jumps of the DG piecewise polynomial basis, e.g. at $x_j$ 
\be
 \widehat u=\frac{u_j^- + u_j^+}{2}, \qquad   {\widetilde u}_x =\frac{{u_x}_j^- + {u_x}_j^+}{2} -\frac{\gamma \; (u_j^+ - u_j^-)}{\Delta x_j}, 
 \label{dgterms}
 \ee
 where $\Delta x_j$ is the length the interval $[x_j, x_{j+1}]$  and 
the rest of the variables are evaluated using their respective 
left and right values on the interval $[x_j, x_{j+1}]$. The constant 
$\gamma$ was set experimentally to $\gamma =200 (p+1)^2$, 
see \cite{Chen}, where 
$p$ is the degree of the polynomial used in the DG approximation.

For the intervals that have a vertex as one of their endpoints, $\widehat u $ is computed as an arithmetic average over all edge values that share the common vertex. Note, that the term containing $(u-\widehat u)$ is an additional a penalty term to enforce the continuity condition at the interval interfaces (and may omitted), but it is the only penalty term that is used at the vertices to enforce the continuity of the solution there. Without it, the numerical solution is 
convergent, but the limiting function is discontinuous at the vertices, 
and thus it is not a solution of the original Poisson problem with 
the stated vertex conditions (continuity plus Kirchhoff constraints).

The flux value $u_x^{(j)}$ at the $j$th vertex is computed from Kirchhoff's condition, 
\be
u_x^{(j)}=-\sum\limits_{\text{all}\; m \neq j} u_x^{(m)}\;.
\ee
Therefore, DG allows one to apply exact Kirchhoff flux conditions in 
the problem formulation whereas the continuity condition is 
enforced via the penalty term.

Finally, for each arc and each interval on these arcs 
$\displaystyle [x_j,x_{j+1}]$ of size 
$\displaystyle \Delta x_j=x_{j+1}-x_j$, the numerical solution 
$u(x)$ is represented as
\be
u(x)=\sum\limits_{m=0}^p c_m \tilde P_m (x),  
\label{dgsol}
\ee
where $\displaystyle \tilde P_m(x)=P_m \Big (-1+ 2 (x-x_j)/ \Delta x_j\Big)$, 
where $P_m(x)$ are standard Legendre polynomials orthogonal on the
interval $x \in [-1,1]$, 
\be
\int\limits_{x_j}^{x_{j+1}} \tilde P_m(x) \tilde P_k(x) \; {\rm d}x = \frac{\Delta x_j}{(2m+1)} \; \delta_{mk}.
\ee
The degree $p$ of the polynomials on each arc and each interval within it 
may vary arbitrarily as communication between neighbouring polynomials 
is achieved via interface or vertex values that can be easily evaluated 
using the above formulas regardless of the polynomial degrees involved.

On the interval $[x_j, x_{j+1}]$ we also use the following properties 
to compute the interface values of the Legendre polynomials, 
\be
 \tilde P_m(x_{j+1})=1, \;\;\; \tilde P_m(x_j)=(-1)^m, \\
\ee
\be
 \tilde P'_m(x_{j+1})=\frac{m(m+1)}{\Delta x_j},\;\;\;  \;\tilde P'_m(x_j)=(-1)^{m+1} \; \tilde  P'_m(x_{j+1}),
\ee
as well as precompute each inner product on the standard $[-1,1]$ interval using identity relating it to the arbitrary interval $[x_j,x_{j+1}]$,
 \be
 \int\limits_{x_j}^{x_{j+1}} \tilde P_m'(x)\tilde P_l'(x) {\rm d}x = \frac{2}{\Delta x_j} \int\limits_{-1}^{1}  P_m'(x) P_l'(x)\;  {\rm d}x.
 \ee
The test function for the DG method on each interval $[x_j,x_{j+1}]$ are $\tilde P_l(x)$,  $l=0,1,2,\dots ,p.$

In the following section we compare the theoretical and numerical error 
estimates for the FD and DG methods.

\section{Numerical results}

Here we validate an implementation of the numerical methods
described in the previous sections for  
the eigenvalue problem, the steady-Poisson equation and the 
time-dependent wave equation. We choose three model graphs:
a {\it pumpkin} graph \cite{berkolaiko17}, a graph obtained from
an electrical grid model, and a graph coming from a laser
based on a random network of optical fibres.

\begin{figure}[H]
\centerline{ \epsfig{file=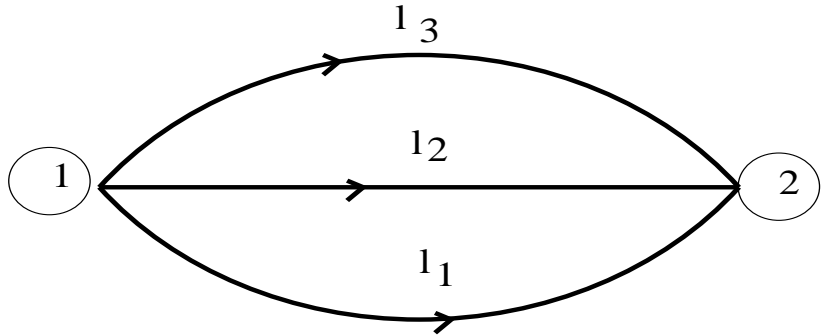,height=5cm,width=12cm,angle=0} }
\caption{A three arc pumpkin graph.}
\label{gg2}
\end{figure}

\begin{figure}[H]
\centerline{
\epsfig{file=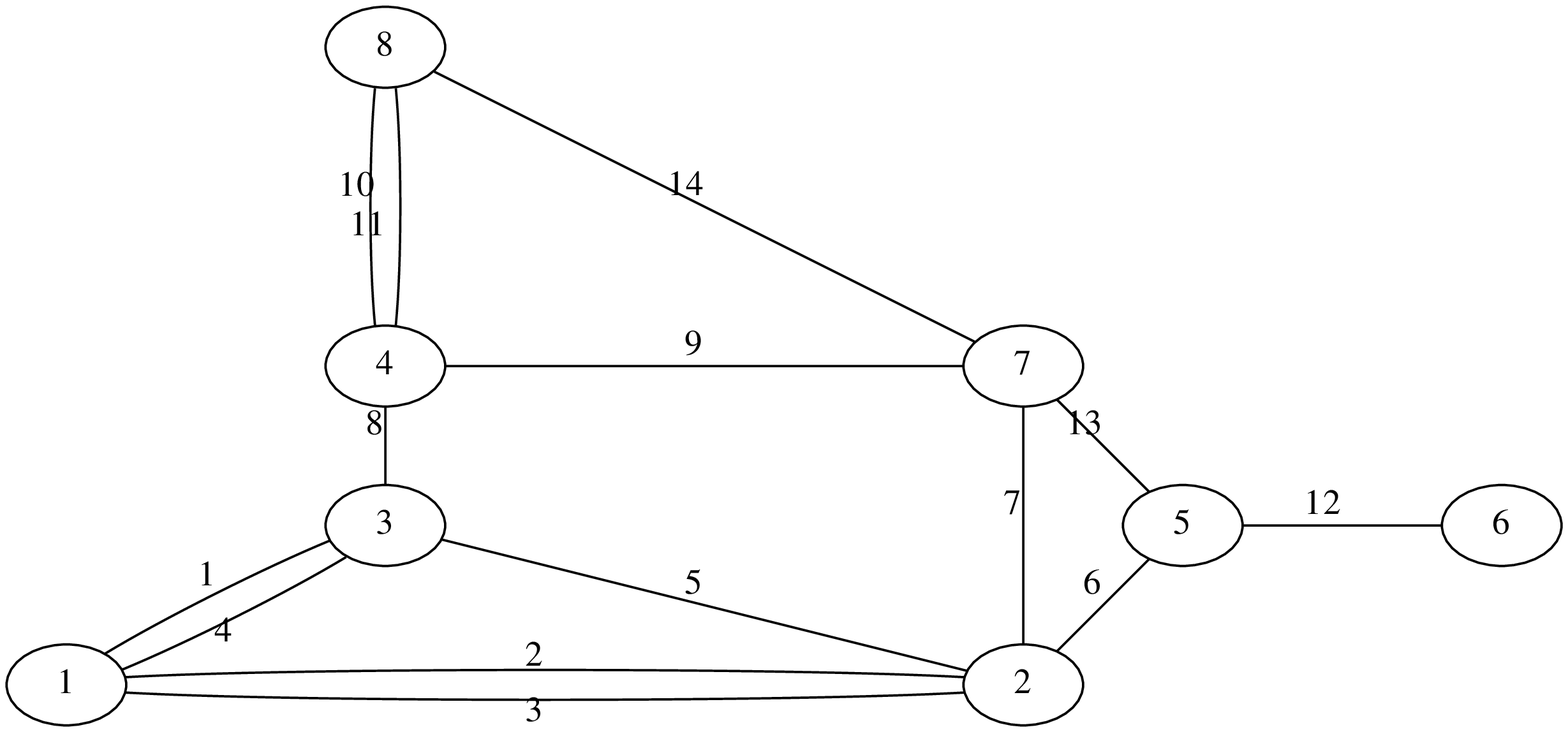,height=5cm,width=10cm,angle=0}
}
\caption{The 14-vertex graph}
\label{gg14}
\end{figure}

\begin{figure}[H]
\centerline{
\epsfig{file=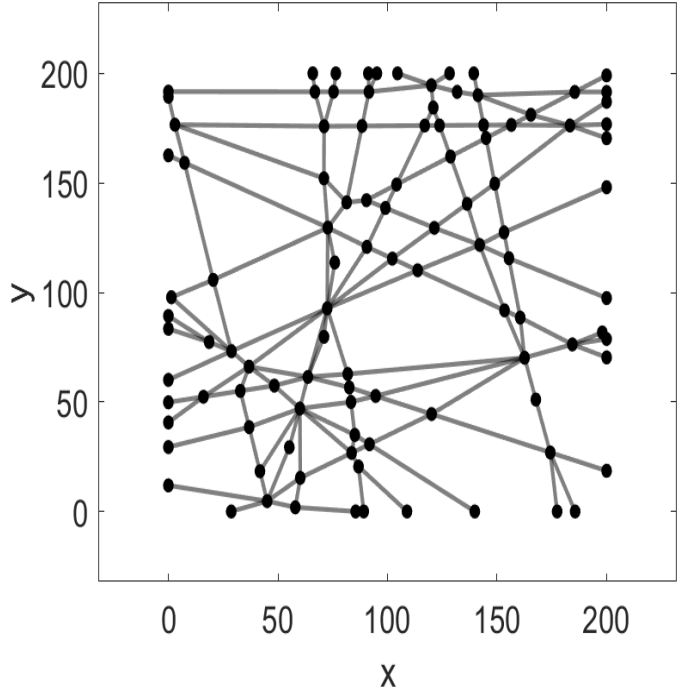,height=10cm,width=10cm,angle=0}
}
\caption{Buffon's needle graph.}
\label{buf1}
\end{figure}

The first graph is a simple three arc pumpkin graph shown in Fig.~\ref{gg2} with arc lengths $l_1=\sqrt{2}, ~ l_2=\sqrt{3},  ~~l_3=\sqrt{5}$, 
A second graph is a 14-vertex graph G14 adapted from the IEEE benchmark test, 
see \cite{g14}. It is shown in Fig. \ref{gg14}. 
The lengths $l_i$, $i=1,\dots ,14$ are given in the table below.
\begin{table} [H]
\centering
\begin{tabular}{|l|c|c|c|r|}
\hline
$l_1$       & $l_2$      &$l_3$& $l_4$       & $l_5$       \\
11.91371443 & 7.08276253 &   6 & 2.236067977 & 4.123105626 \\ \hline
$l_6$       & $l_7$      &$l_8$& $l_9$       &  $l_{10}$   \\
1.414213562 &  2         & 1   & 4.7169892   & 4.472135955 \\ \hline
$l_{11}$    &  $l_{12}$  & $l_{13}$ & $l_{14}$ &   \\
2           & 2          & 1.414213562 & 4.472135955 & \\
\hline
\end{tabular}
\caption{The lengths $l_i$ for the graph G14. }
\label{tab2}
\end{table} 
We also consider an example of a random graph, a Buffon's needle graph with 165 arcs and 104 vertices, see Fig. \ref{buf1}. We generated a version of Buffon's needle graph  using a code developed by Michele Gaio, \cite{Gaio}. To create the graph, $n$ random points are uniformly generated 
on a square with diagonal $L$. A straight line segment (needle) with 
random angle and length is drawn from each point. The needles are 
all given the same diameter $D$. The intersections of line segments 
become the vertices of the graph. If two vertices are within a distance 
$r=\left(\frac{12D}{L}\right)^2$ of each other, they are combined 
into a single vertex. The arcs of the graph are the lengths along 
the needles connecting the vertices. The resulting graph appears as 
a series of scattered needles. By construction, most vertices are degree four, with on average 15\% of  vertices of degree six.

\subsection {Eigenvalues and Weyl's law}

To illustrate the spectral algorithm for finding eigenvalues using exact eigenvectors given by Fourier modes of equation  (\ref{simple}),
we outline the steps involved on a pumpkin graph  (Fig. \ref{gg2}). 
Enforcing continuity and Kirhchoff's law at vertices $1$ and $2$, respectively, 
gives
\begin{eqnarray}
B_1=B_2=B_3,\\
A_1 + A_2 + A_3 = 0,\\
A_1 s_1 + B_1 c_1 = A_2 s_2 + B_2 c_2 = A_3 s_3 + B_3 c_3,\\
A_1 c_1 - B_1 s_1 + A_2 c_2 - B_2 s_2 + A_3 c_3 - B_3 s_3=0,
\end{eqnarray}
where $s_1= \sin k l_1, c_1=\cos k l_1$, etc. 
This yields the following linear system
{\small
\be\label{ling2}
M(k) \begin{pmatrix} A_1 \cr B_1 \cr A_2 \cr B_2 \cr A_3 \cr B_3 \end{pmatrix}
=
\begin{pmatrix} 0 \cr 0 \cr 0 \cr 0 \cr 0 \cr 0 \end{pmatrix} .
\ee
where the matrix $M(k)$ is
\be\label{mk}
M(k) \equiv \begin{pmatrix}
.  & 1 & . & -1 & . & . \cr
.  & 1 & . &  . & . & -1 \cr
1  & . & 1 &  . & 1 & . \cr
s_1  & c_1 & -s_2 &  -c_2 & . & . \cr
s_1  & c_1 & .    &  .    & -s_3 & -c_3 \cr
c_1  & -s_1 & c_2 &  -s_2 & c_3 & -s_3
\end{pmatrix}
\ee
}
Setting the inverse condition number of the matrix $M(k)$ to $0$, we obtain the
equations for the resonant frequencies $k$, and then compute the eigenvalues using relation  $\lambda=-k^2$. The first three nonzero eigenvalues  and corresponding eigenvectors
$V^q =(v^q_1,v^q_2,v^q_3)^T$ are given in the Table \ref{tab1}.
\begin{table} [H]
\centering
\begin{tabular}{|l|c|c|c|c|}
\hline
$q$ & $-k_q^2$    & $A^q_1 ~/~B^q_1$ & $A^q_2 ~/~B^q_2$ & $A^q_3 ~/~B^q_3$\\ \hline
1 & -2.395998     &  -0.24204 &-0.53262 & 0.77466\\
  &      &  -0.12486 &-0.12486 & -0.12486\\ \hline
2 & -3.057162     &0.20191    &0.03291  & -0.23481\\
  &        &-0.58105   &-0.58105 & -0.58105\\ \hline
3 & -4.067077     &0.85001    & -0.69799  & -0.15202\\
  &      &0.123927   & 0.123927  & 0.123927\\
\hline
\end{tabular}
\caption{First three nonzero eigenvalues and eigenvector components for the pumpkin graph.}  
\label{tab1}
\end{table}

\begin{figure}[H]
\centerline{
\epsfig{file=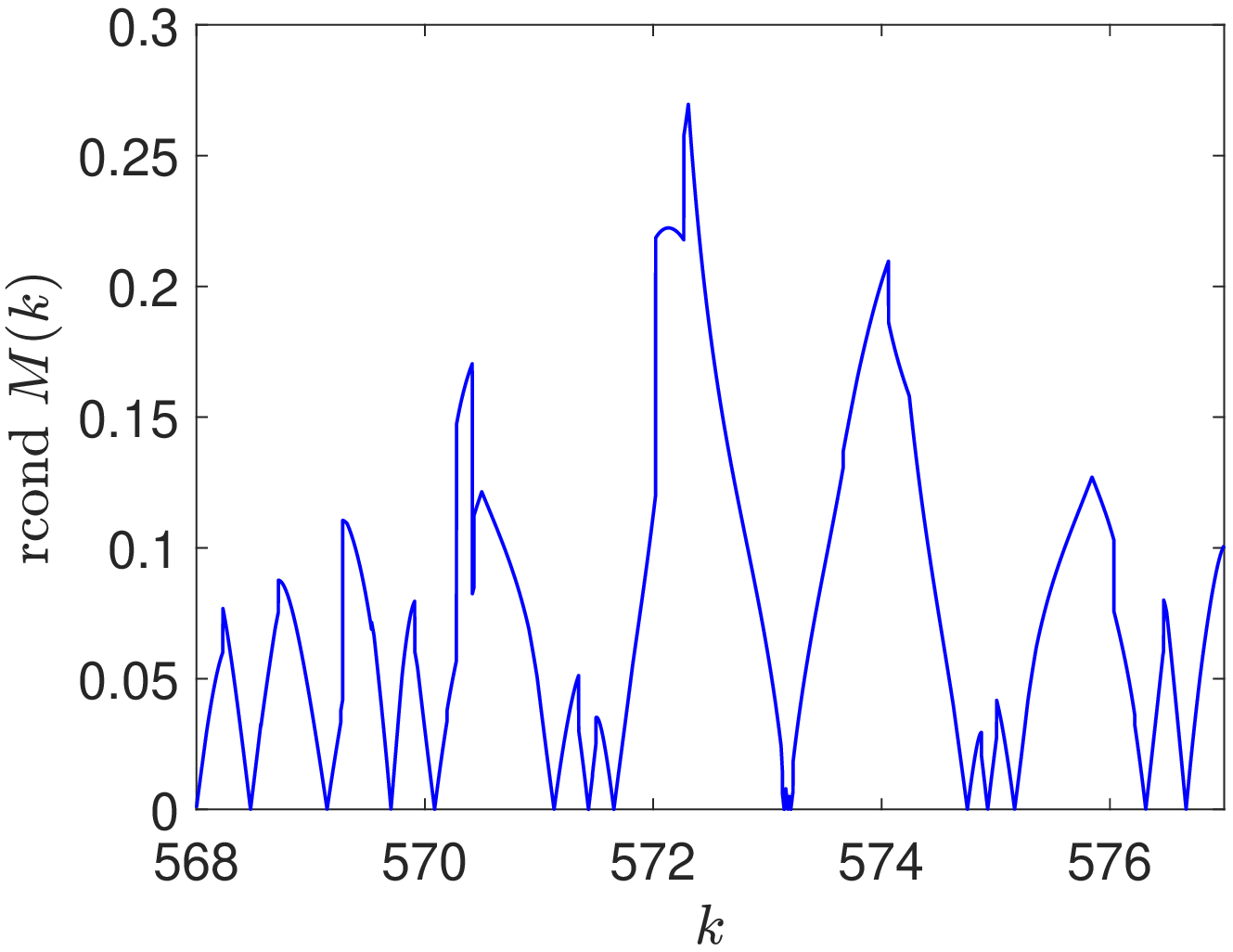,height=7 cm,width=6 cm,angle=0}
\epsfig{file=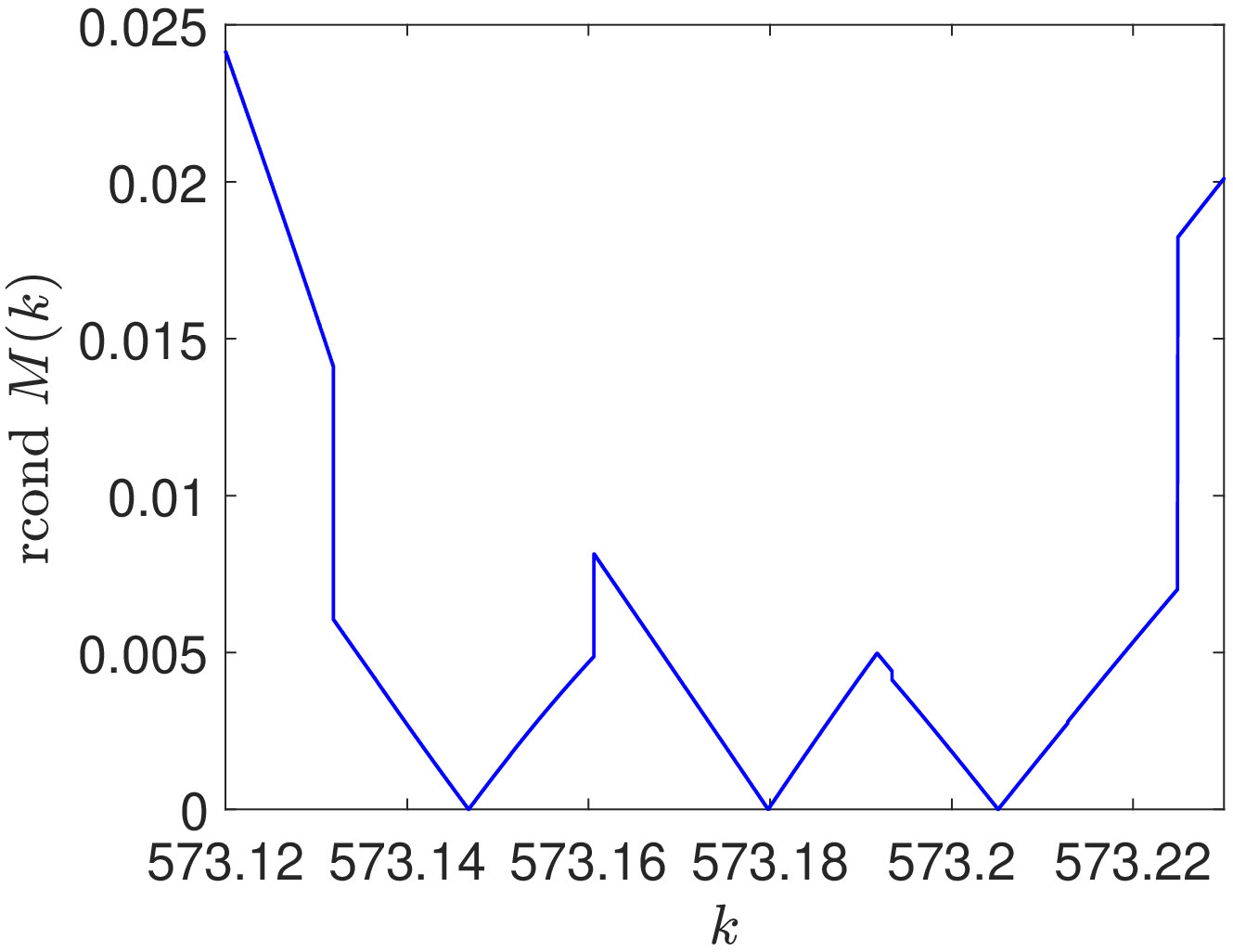,height=7 cm,width=6 cm,angle=0}
}
\caption{Pumpkin graph: plot of rcond(k) for $k \in [568,577]$ (left panel),
blow-up of the plot (right panel).}
\label{rcond2}
\end{figure}

To compute the resonant frequencies we estimate the minimum spacing graphically and split the range of $k$ into intervals of about $1/10$ of the estimated 
minimum spacing. Each minimum of $r(k)$ is found using 
a Brent-like minimisation algorithm that we implemented, e.g. 
a secant method with bracketing see \cite{linesearch} for example.
A more efficient algorithm that we did not pursue here  would be to 
use the sawtooth nature of the $\text{rcond(M(k))}$ functions as shown in 
Fig. \ref{rcond2} and alternate between min and max 
search algorithms to find peak and valley values only.

To illustrate the practicality and robustness of our algorithm, we
show rcond(k) for $k \in [568,577]$ in Fig. \ref{rcond2}. Such
high-order eigenvalues and eigenvectors
can be computed easily with the spectral method as opposed to FD or DG,
for which the spatial structure of the eigenvectors need to be resolved.
The spacing between the resonant $k$s is fairly regular except
at exceptional locations, such as the one shown in Fig. \ref{rcond2}
where the $k$ in $[572,574]$ need a finer bracketing, see right panel
of the figure. Nevertheless, the line minimisation yields the
estimates of $k$ at machine precision
$$573.14678431204834,~~ 573.17977474169390,~~ 573.20510976082187~ .$$

{\bf Weyl's law}

\begin{figure}[H]
\centerline{
\epsfig{file=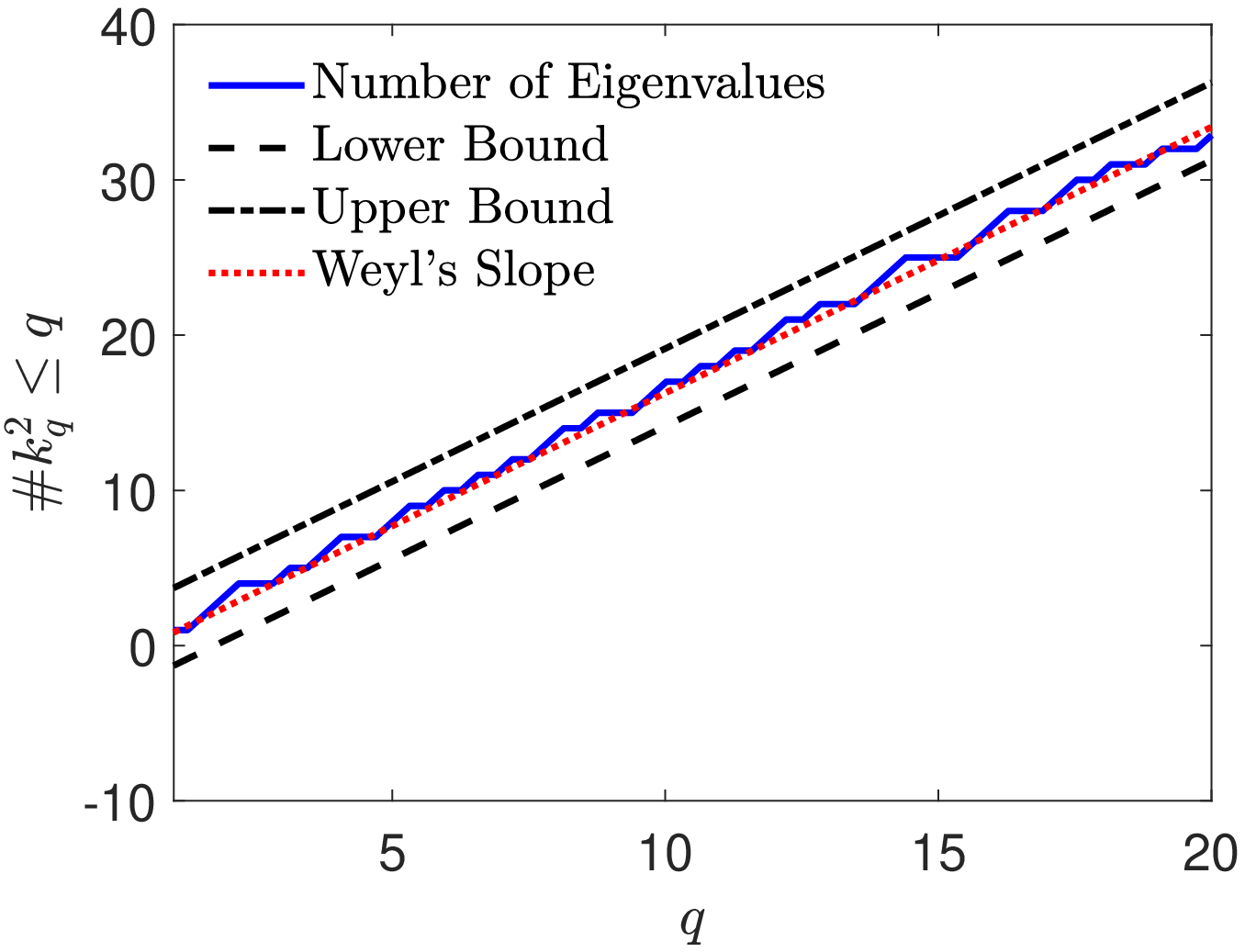,height=5cm,width=4 cm,angle=0}
\epsfig{file=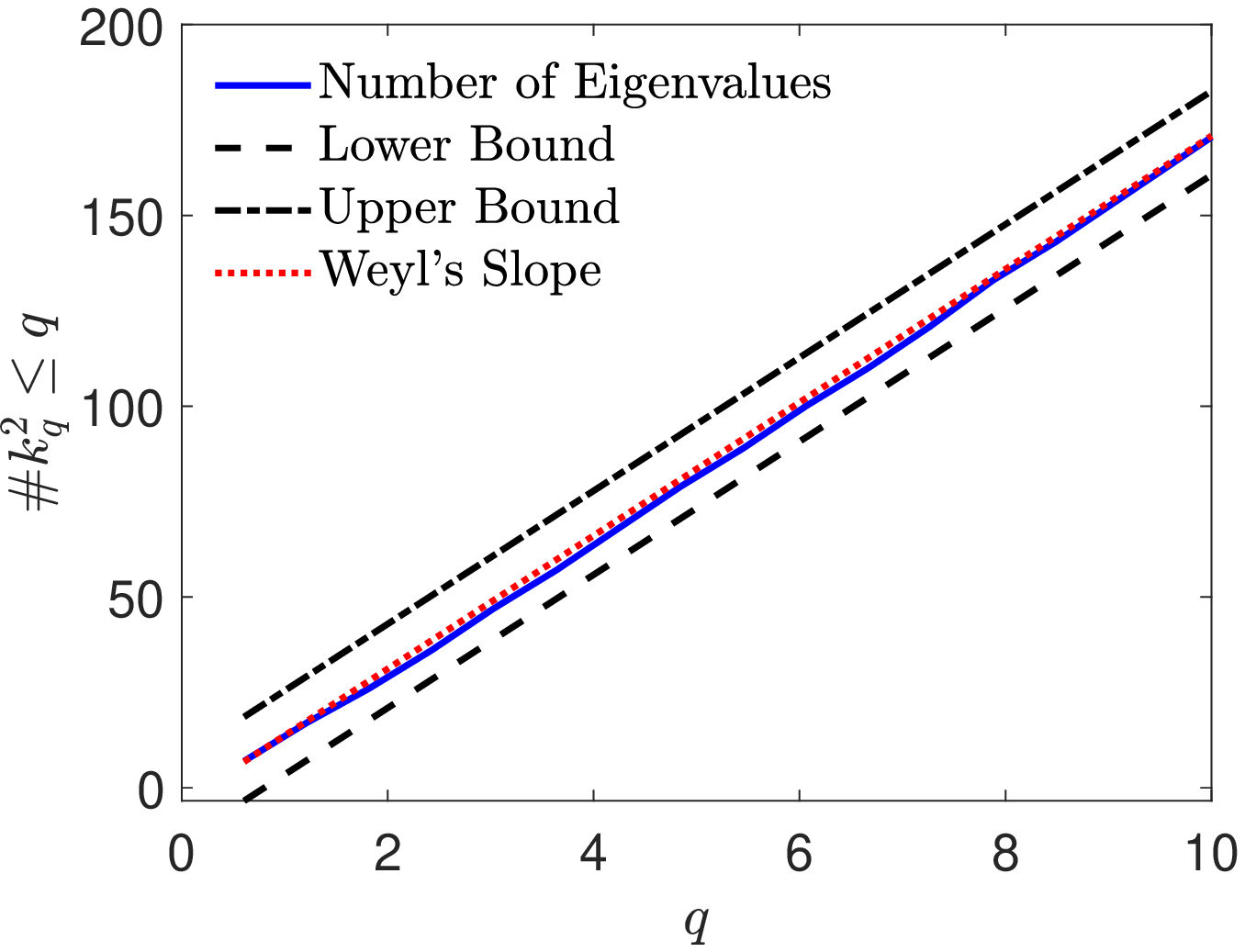,height=5cm,width=4 cm,angle=0}
\epsfig{file=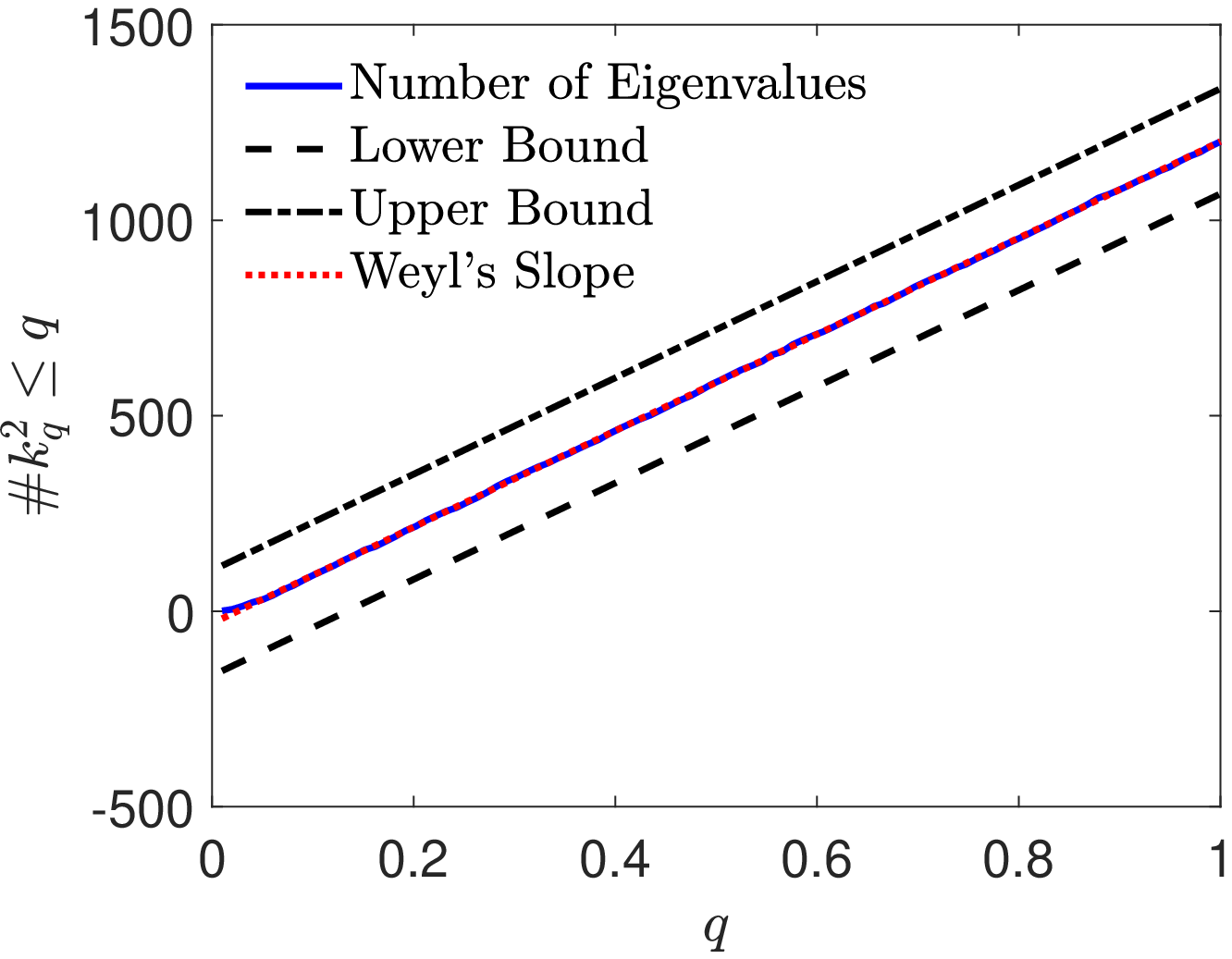,height=5cm,width=4 cm,angle=0}
}
\caption{Distribution of the resonant frequencies vs Weyl' law 
estimates for the Pumpkin (left),  G14 (middle)  and Buffon's (right) graphs. }
\label{weylaw}
\end{figure}

The eigenvalues $k^2$ of the generalised Laplacian (\ref{helm}) of
graphs follow Weyl's law \cite{berkolaiko17}. For a positive real number $q$, 
the number of eigenvalues $k^2 \le q^2$ is given by Weyl's estimate
\be\label{weyl} q L / \pi , \ee 
where $L$ is the sum of the lengths of all the arcs of the graph. 
In addition, the distribution of eigenvalues is also bounded from above 
and below by two lines with the same slope
\be\label{weylbounds} q L / \pi - \left| E \right| \leq \;  \# \{ k \big| k^2 \leq q^2 \} \; \leq q L / \pi + \left| V \right|,   \ee
 as illustrated in Fig. \ref{weylaw}. Here, $\left| E \right|$ represents the number of arcs, and $\left| V \right|$ represents the number of vertices. 
As seen in  Fig. \ref{weylaw}, the distribution of eigenvalues follows Weyl's law with $L=5$ (pumpkin graph), 40 (G14) and 2700 (Buffon).

\begin{figure}[H]
\includegraphics[scale=0.5]{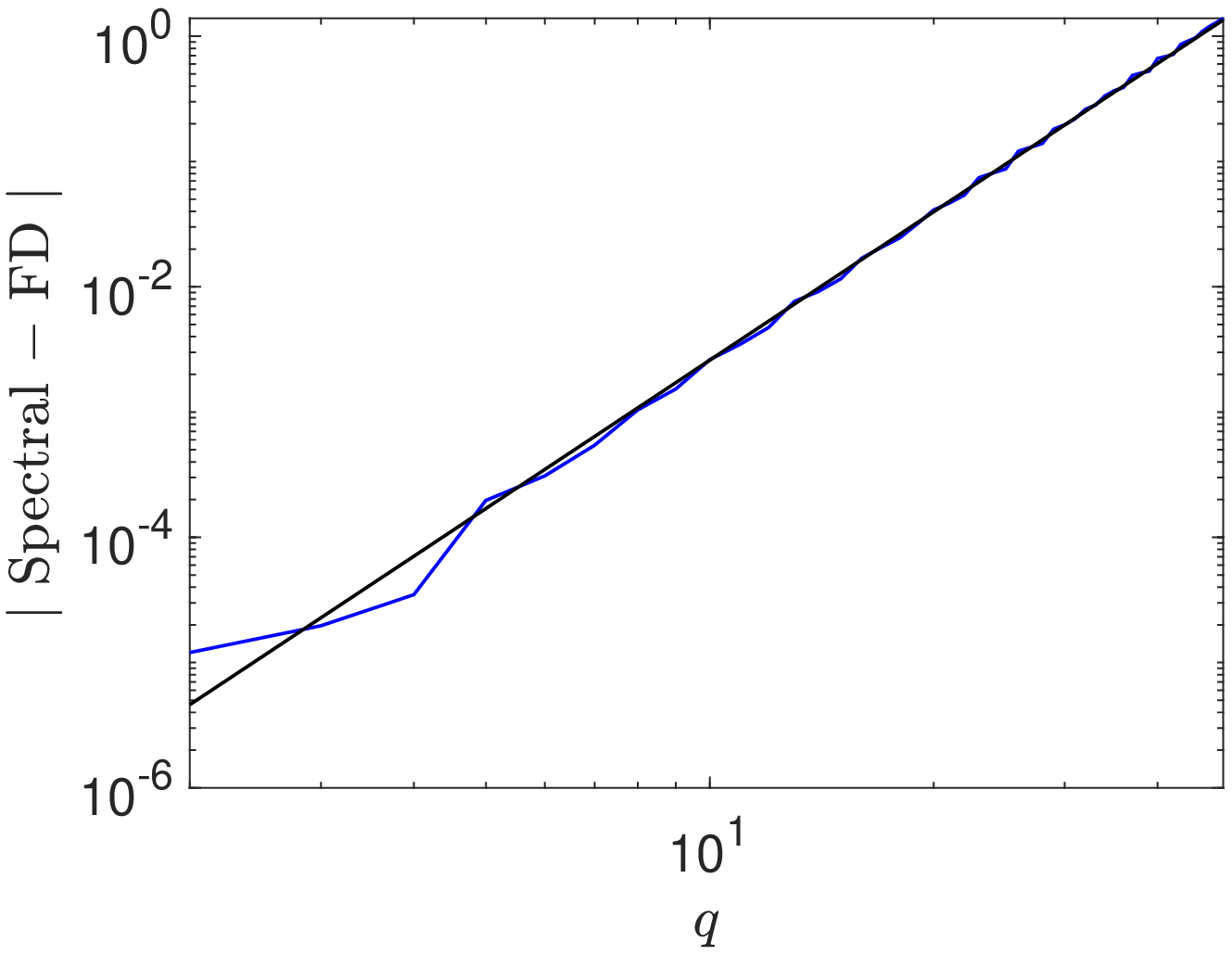}
\includegraphics[scale=0.5]{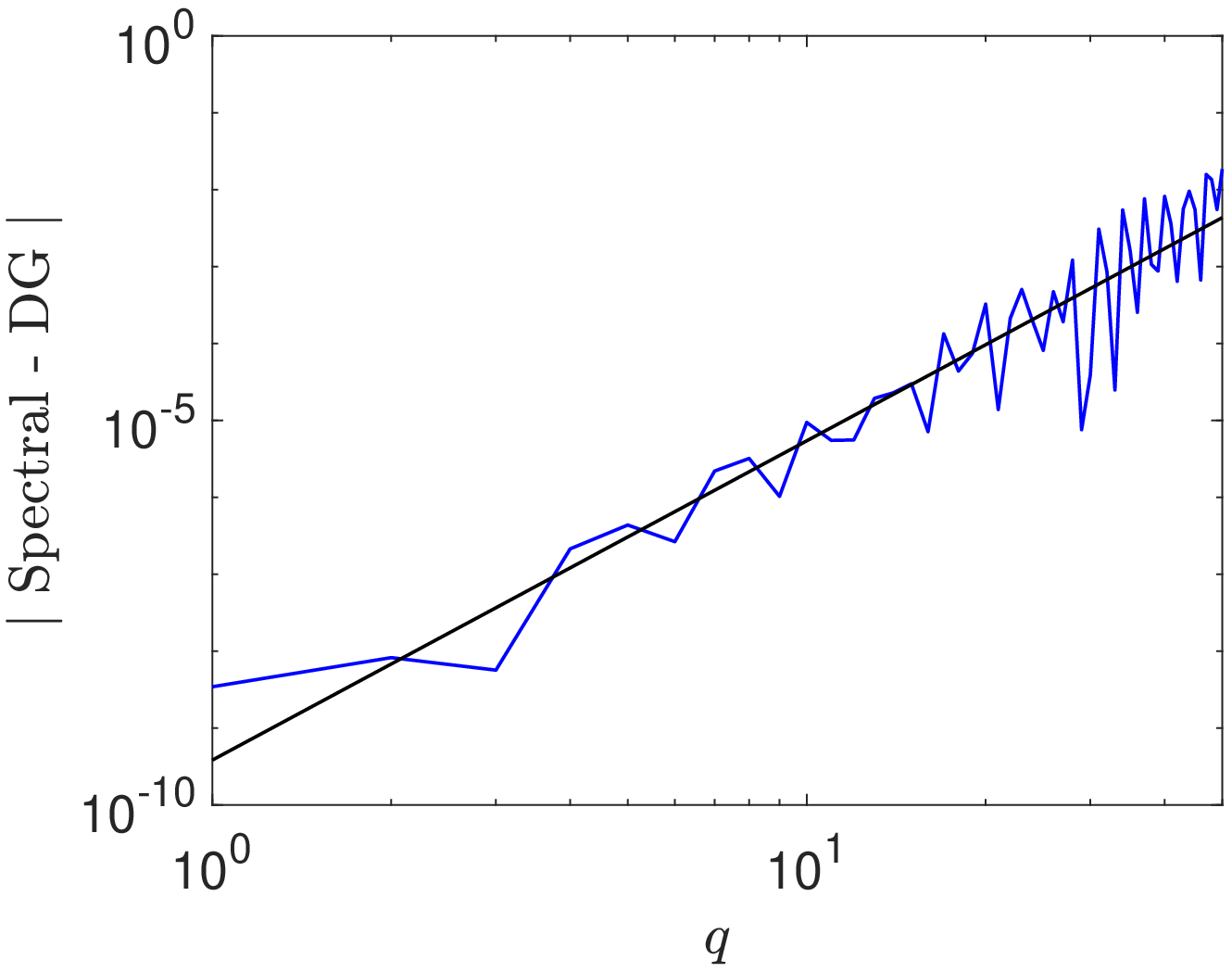}
\caption{Absolute errors for the $q$th eigenvalue vs $q$ 
for the FD (left panel) and
the DG (right panel) algorithms.}
\label{eigs_compare}
\end{figure}

Since root-finding accuracy is independent of the $k$-value 
for the spectral algorithm, it is superior to both FD and 
DG methods, as their accuracy depends on spatially resolving 
the eigenvector. For example, on the pumpkin graph, 
for $\displaystyle k=573.18$ shown in Fig. \ref{rcond2}, 
one has to resolve the wavelength $10^{-3}$.
This would require a spatial step $\Delta x \approx 10^{-12}$ 
for the error to reach machine precision, see equation \ref{errEigFD}
below.
For DG with 5th-order polynomial, the sub-intervals
should be about $h= 10^{-3}$ for the error to reach machine precision.
For fixed mesh size, the consecutive eigenvalues become less and less precise
for both FD and DG methods, as shown in Fig. \ref{eigs_compare} for the first 
30 resonant frequencies for the pumpkin graph. 
Fig. \ref{eigs_compare} shows the errors in the
eigenvalues for the FD and DG methods, using the spectral
estimation as an exact value. 
Note the linear scalings of the absolute errors. For the FD, we have
\be \label{errEigFD} \log_{10}({\rm abs. ~ error}) = 3.93  \log_{10}(q)   -  6.52 \ee
The eigenvalues of the tridiagonal matrix that approximates
the second derivative 
are  $$-(4/\Delta x^2) \sin (\pi k \Delta x/2)^2 \approx - k^2 \pi^2 + O( k^4 \Delta x^2) $$  
in a Taylor expansion around $\Delta x=0$.
This explains the scaling $k^4$ observed for the error of the FD
method.

For the DG, the equation of the line is
\be \log_{10}({\rm abs. ~ error})=4.15  \log_{10}(q) - 9.42 ~ . \ee

\begin{figure}[H]
\centerline{ 
\epsfig{file=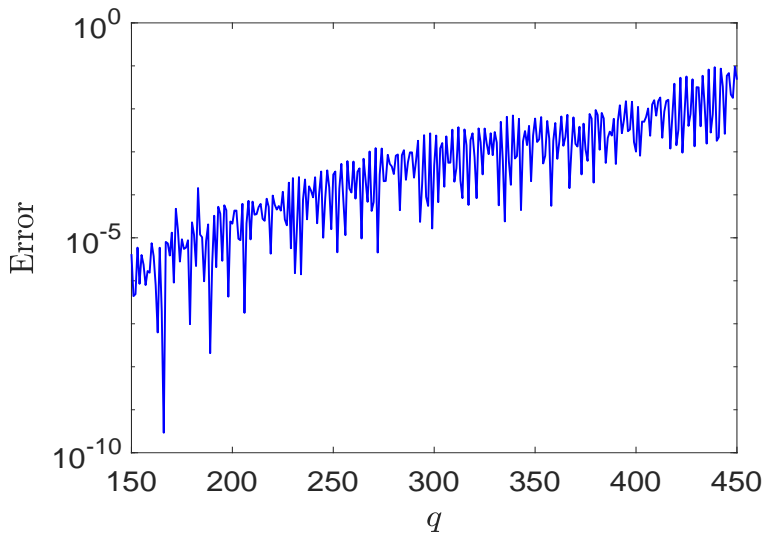,height=7 cm,width=10 cm,angle=0}
}
\caption{Pumpkin graph: error in DG estimation of eigenvalues of
high order, $h=0.01,~~ p=5$. }
\label{dg2}
\end{figure}

The DG method can be used to estimate the eigenvalues of high-order for
the generalised Laplacian (\ref{helm}). Fig. \ref{dg2} shows the
error as a function of $q$.
As expected, accuracy decreases as the order of the eigenvalue increases.
For $h=0.01$,  the error for $q=150$ is $10^{-5}$, and the error reaches 
$10^{-1}$ for $q=450$.

Since the form of the eigenvectors of the generalised Laplacian (\ref{helm})
are known exactly, the only approximation required is to find 
the real zeros of a scalar real function. This can be done 
for thousands of eigenvalues within machine precision. 
In contrast, for the FD and DG methods, computing the eigenvalues 
accurately requires the accurate resolution of the corresponding 
eigenvectors. This is much more computationally expensive, even 
with the DG method that allows an arbitrary polynomial order of 
accuracy, $O(h^{p+1})$, where $h$ is the interval size 
in splitting the arcs and $p$ is the smallest polynomial 
order used for approximating the solution on each interval.

In the next two sections,  we discuss solutions of the Poisson 
and wave equations using the three methods under consideration: 
spectral Fourier, FD, and DG methods. 

\subsection {Poisson equation}

Consider Poisson's equation 
\be\label{poisson} 
\hat \Delta U = F 
\ee
where $\hat \Delta$ stands for the generalised Laplacian on the graph (\ref{helm}).

Using the spectral decomposition, expand both $U$ and $F$ in terms of the 
eigenvectors $V^q$ of the equation (\ref{helm}), 
$\displaystyle \Delta V^q =- k_q^2 V^q$, 
\be
U = \sum_{q=0}^{\infty} \alpha_q V^q, ~~~~F=\sum_{q=0}^{\infty}\beta_q V^q.
\label{expan2}
\ee
Since  0 is an eigenvalue
of the Laplacian with corresponding constant eigenvector, the 
compatibility condition (Fredholm alternative) requires $\beta_0=0$. 
The rest of the unknown coefficients $\alpha_q$ are obtained by
projecting onto the eigenvectors $V^q$ using the scalar product
(\ref{vivi})
\be\label{sollap}
\alpha_q = -{\beta_q \over k_q^2}, \;\; q=1,2,3, \dots \ee

Since the eigenfunctions are orthogonal, Bessel's identity allow one to measure the error due to series truncation in terms of the decay of the expansion coefficients. 
  
\begin{figure}[H]
\centerline{
\epsfig{file=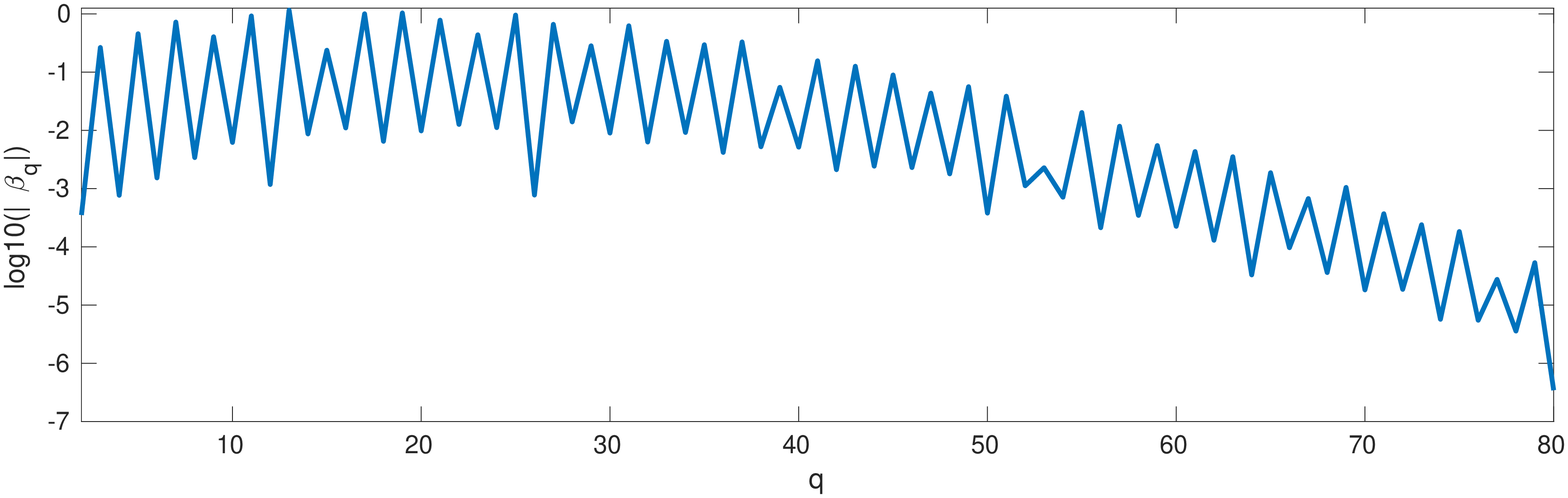,height=6cm,width=7 cm,angle=0}
\epsfig{file=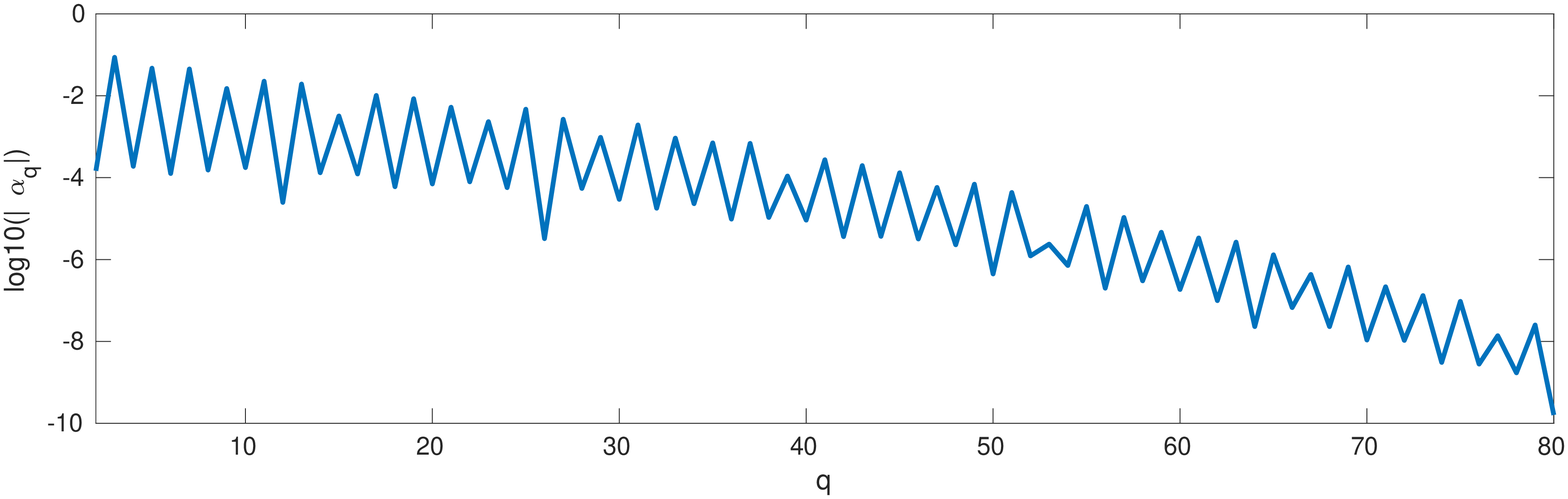,height=6cm,width=7 cm,angle=0}
}
\caption{Base 10 logarithms of the amplitudes $\beta_q$ (left) and
$\alpha_q$ (right) for the function $g$ given as (\ref{dgauss}) on arc 2
of the pumpkin graph.} 
\label{pumkGauss}
\end{figure}

In Fig. \ref{pumkGauss} we show the accuracy of the solution for the pumpkin graph in terms of the accuracy of the function $f$ expansion,   $\log10(\beta_q)$ (left)  and 
$\log10(\alpha_q)$ (right), respectively, in terms of the number of terms in partial sums in (\ref{expan2}). The initial condition is a derivative of 
the Gaussian on arc 2 (that has zero average) and zero on the other arcs.
\be\label{dgauss}
g(x) =  -2 {x-x_0 \over  s^2} \exp(-{(x-x_0)^2 \over  s^2}), \ee
 with $x_0= 0.865  ~~s=0.15$.

The problem has an exact solution, up to an arbitrary constant
\begin{align}  \label{exactPoiss}
& u_1(x) = -0.04431134627263788525 + 0.06267946415920350157 x, \\
& u_2(x) = 0.08862269254527577050 - 0.10232143801160161859 x + {3 \sqrt{ \pi} \over 40} 
{\rm erf}(  {200 x- 173 \over 30} ) , \\
& u_3(x) = -0.04431134627263788525 + 0.03964197385239476854 x  .
\end{align}

\begin{figure}[H]
\centerline{
\epsfig{file=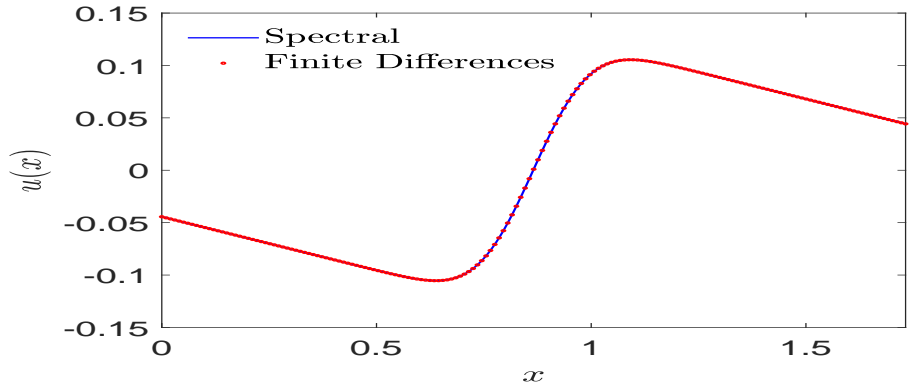,height=5cm,width=12 cm,angle=0}
}
\caption{
Solutions of Poisson's equation (\ref{poisson}) with a right hand side
given by (\ref{dgauss}) in arc 2 and zero in the other arcs. }
\label{edge2_solution}
\end{figure}

The consistency of the three numerical methods with the exact solution 
is shown in Fig. \ref{edge2_solution}. The arbitrary constant
is determined by a least square fit of the difference between the exact
and the numerical solution.

Here we solve (\ref{poisson}) for the pumpkin graph with 
$\Delta x = 0.01$ and  a right hand side $F$ given by (\ref{dgauss})
placed on arc 2 with length = $\sqrt{3}$ and $0$ on the other
edges. The parameters are $x_0=0.865$ and $s=0.15$.

The convergence of the spectral amplitudes $\alpha_q$ and $\beta_q$ 
vs the number of modes $q$ is shown in Fig. \ref{pumkGauss}.  
The semilog plot shows that the coefficients $\beta_q$ and
$\alpha_q$ decay exponentially with $q$ (spectral accuracy),
as expected for the Fourier expansion of analytical functions. For example, 
an expansion with fifty modes allows one to reach a solution with 
accuracy (absolute error) of $10^{-8}$.
\begin{figure}[H]
\centerline{
\epsfig{file=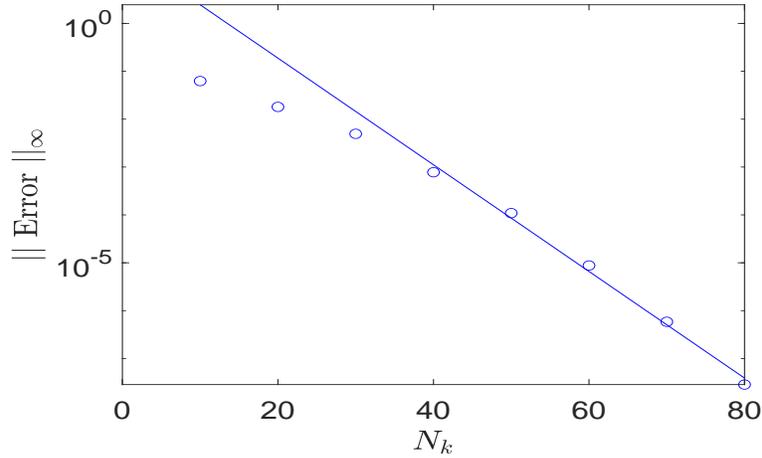,height=6cm,width=10 cm,angle=0}
}
\caption{$L_\infty$ norm of the difference between the spectral solution
and the exact solution (\ref{exactPoiss}) as a function of the number of 
modes $N_k$ used in the expansion.}
\label{errorSpec}
\end{figure}
To confirm these results, we plot in Fig. \ref{errorSpec} in log-linear
scale the $L_\infty$ norm of the difference between the spectral solution 
(\ref{sollap})
and the exact solution  (\ref{exactPoiss}) as a function of $N_k$ the
number of modes used in the expansion. As expected, we observe
an exponential decay of the error.

For both the FD and DG methods, solving the 
Poisson equation (\ref{poisson}) reduces to solving the linear system $\displaystyle A \tilde u =b$,where $A$ represent the discrete Laplacian approximated using the FD or DG method and where $b$ represents the 
strong and weak approximations of the function $F$, respectively, 
as described in the previous section. The $\tilde u$ represents the unknown discrete values of $u$ or the expansion coefficients in Legendre polynomials 
for the DG method.

The matrix $A$ is singular due to the zero eigenvalue corresponding
to a constant eigenvector and a pseudo-inverse is used to find the
solution without the arbitrary constant belonging to the null 
space of A. We rewrite $A$ in its reduced SVD
form,  $A=\hat{U} \hat{\Sigma} V^T$ and solve for $u$,  
\be
u = A^{\dagger} F, \ee
where  $A^{\dagger} = V \hat{\Sigma}^{-1} \hat{U}^T$.

\begin{figure}[H]
\centerline{
\epsfig{file=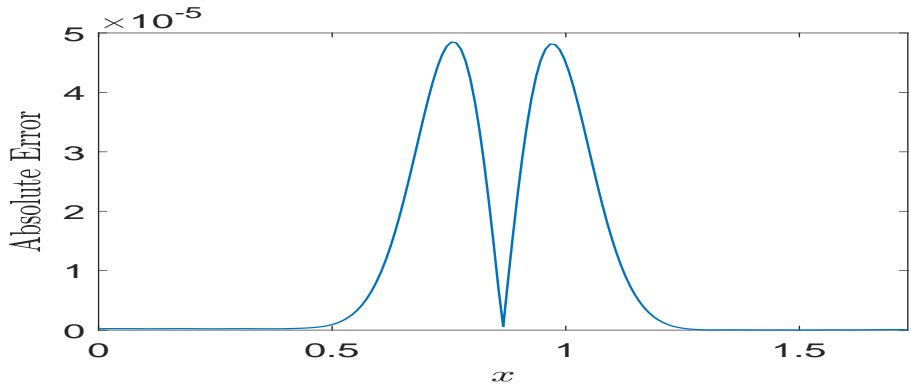,height=5cm,width=12 cm,angle=0}
}
\caption{ 
Errors of the FD method for the 
solution of Poisson's equation (\ref{poisson}) with a right hand side
given by (\ref{dgauss}) in arc 2 and zero in the other arcs.}
\label{erDGFDpoisson}
\end{figure}

The error for the FD method is shown in Fig.  \ref{erDGFDpoisson}.

\begin{figure}[H]
\centerline{
\epsfig{file=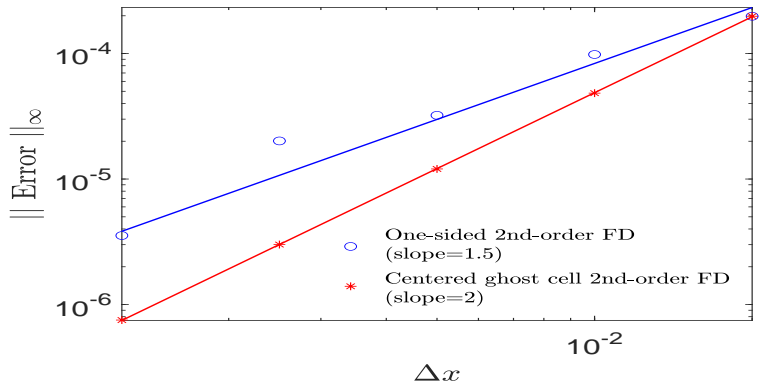,height=5cm,width=10 cm,angle=0}
}
\caption{
Plot of the error vs. $\Delta x$ for two approximations of the
Kirchoff vertex condition:  centered FD with ghost points 2nd order 
(red online) and one-sided FD approximation (blue online)}
\label{DF21}
\end{figure}

To conclude on the FD, we show in Fig. \ref{DF21} the error for
the standard vertex approximation used by \cite{Arioli}
and the method we used, i.e. centered FD with ghost points. Clearly
the former is first order while the latter is second order.

The DG approximation converges at the optimal convergence rate $O(h^{p+1})$ as illustrated in table below that provides the errors in 
$L_{\infty}$-norm. The error is maximal near the inflection point of 
the solution on the second edge. It is lower by an order 
of magnitude or two away from the inflection point and on the other 
two edges. The condition number of 
the matrix grows exponentially with the size $n$ of the matrix 
and the error saturates 
for $p=5$ and $h=0.01$, see Table \ref{DGerror}.
At this point, the SVD loses precision,  
as confirmed by the exact computation in rational arithmetic using 
Mathematica, which gives an error for 
$p=5, h=0.01$ of $4.9~10^{-13}$ instead of $3.3~10^{-10}$.
To avoid this loss of accuracy in finite precision, 
one could use a preconditioner as suggested by \cite{Arioli}.
\begin{table}[H]
\hspace*{-0.2cm}
    \begin{subtable}[h]{0.45\textwidth}
        \centering
       \begin{tabular}{|l|c|c|c|c|}
       \hline
       \diagbox[innerleftsep=.2cm,innerrightsep=2pt]{p}{$h$}
         & 0.1    &  0.01\\ \hline
       1 & $7.0~10^{-3}$  &  $6.0~10^{-5}$\\
         \hline
       2 & $6.9~10^{-4}$     & $7.0~10^{-7}$   \\
        \hline
       3 & $5.7~10^{-5}$     & $5.0~10^{-9}$  \\
        \hline
        4 & $6.0~10^{-6}$     & $8.5~10^{-11}$  \\
         \hline
        5 & $3.3~10^{-7}$     & $3.3~10^{-10}$  \\
         \hline 
       \end{tabular}
       \caption{ $O(h^{p+1})$ numerical error  }
       \label{tab:week1}
    \end{subtable}
   \hspace*{-0.4cm}
    \begin{subtable}[h]{0.55\textwidth}
        \centering
       \begin{tabular}{|l|c|c||c|c|}
       \hline
      p   & n    & $cond(A)_2$ & n    & $cond(A)_2$ \\ \hline
            1 & 106  & $1.3 ~10^{5}$  & 1078   & $1.3~10^{7}$   \\
         \hline
       2 & 159  & $5.9~10^{5}$  & 1617  & $6.0~10^{7}$   \\
        \hline
       3 & 212  & $1.0~10^{6}$  & 2156  & $1.0~10^{8}$   \\ 
        \hline
        4 & 265  & $2.4~10^{6}$  & 2695   &$2.5~10^{8}$   \\ 
         \hline
        5 & 318  & $3.5~10^{6}$  & 3234   &$6.6~10^{8}$   \\  
         \hline 
       \end{tabular}
        \caption{ $cond(A)_2$ vs degree $p$ of polynomial}
        \label{tab:week2}
     \end{subtable}
     \caption{DG solution errors in $L_{\infty}$-norm and 
condition number in $L_2$-norm for the matrix $A_{n\times n}$. }
     \label{DGerror}
\end{table}

\subsection {Wave equation}

Here, we consider a telegrapher's equation that generalises both 
the wave and heat equations by adding dispersive and damping terms when $\gamma \ne 0$   and  $\beta \ne 0 $, respectively,  
\be\label{dampweq}
\alpha U_{tt} +\beta U_{t} - {\tilde \Delta} U+ \gamma u=0  . \ee

The eigenvectors $V^q$ of the generalised Laplacian form a 
complete basis of eigenvectors of $L^2(\Omega)$.
It is then natural to expand $U$ as
$$ U= \sum_q \alpha_q(t) V^q .  $$
Substituting the expansion into (\ref{dampweq}) we obtain the following 
evolution equation for the coefficients $\alpha_q$
\be\label{evok}
\alpha {\ddot \alpha_q } + \beta {\dot \alpha_q }
+ k_q^2 \alpha_q + \gamma  \alpha_q =0 . \ee
Given the initial condition this equation can be solved exactly for 
each Fourier amplitude $\alpha(t)$.

Note that when $\alpha=0$, equation (\ref{dampweq})
reduces to the generalised heat equation on the metric graph. 
The evolution of the coefficients $\alpha_q$ is given by
\be\label{evokh}
{\dot \alpha_q } + k_q^2 \alpha_q + \gamma  \alpha_q =0 , \ee
where we chose $\beta=1$ for simplicity.
The solution of (\ref{evokh}) is
\be\label{sevokh}
\alpha_q (t) = \alpha_q (0) \exp\left [ {-( k_q^2 + \gamma)t} \right ]  . \ee
Exactly as for the one-dimensional heat equation, the large $k_q$ coefficients
decay fast and only $\alpha_q$ for the smallest $k_q$ can be observed as time
evolves. It is then particularly important to estimate this first
eigenvalue, the so-called spectral gap \cite{berkolaiko17}. It could explain
the oscillations observed in gas networks, see \cite{cherkov}.

Both FD and DG methods applied directly to (\ref{dampweq}) would require one 
to solve the ODEs for the unknown solution values and expansion 
coefficients $c_m(t)$ in terms of the Legendre polynomials (\ref{dgsol}) 
of the form 
\be
M \ddot C(t) +L\: \dot C(t) +K \: C(t) + \Gamma \: C(t)=0.  
\label{telecoef}
\ee  

\begin{figure}[H]
\centerline{
\epsfig{file=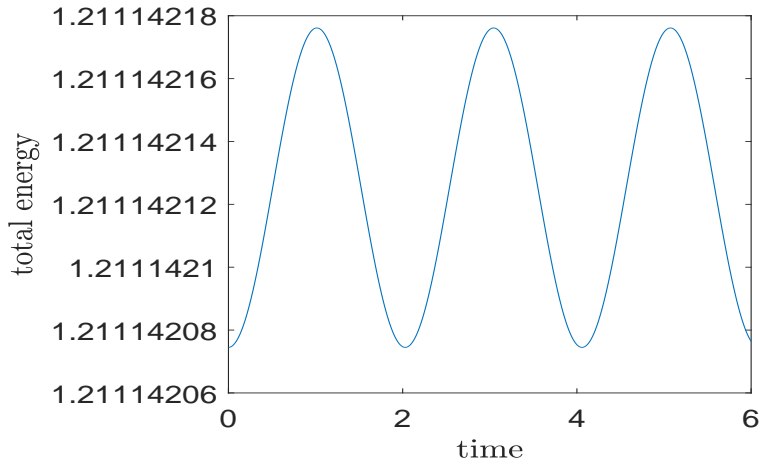,height=6cm,width=10 cm,angle=0}
}
\caption{Pumpkin graph: total energy for the wave equation using the second 
eigenvector $V^2$ ($k_2=1.5479012370538900$) as initial condition. The 
time and space steps are $\Delta t = 5 \times 10^{-4}$ and 
$\Delta x = 6.25 \times 10^{-4}$ . }
\label{energy_mode1}
\end{figure}

Below we illustrate the behavior of the total energy over time for the wave equation on the pumpkin graph.
On a network with $m$ arcs, the energy $E$ is given by
\be
E = \sum_{i=1}^m  \int_{0}^{l_i} \frac{1}{2} \left( \left(\frac{\partial u}{dt}\right)^2 + \left(\frac{\partial u}{dx}\right)^2 \right) \,dx
\ee
for solution $u(x,t)$ and is constant in time. For the pumpkin graph, the 
energy of the second mode is computed to machine precision as 
$E_2 = 1.211142452041264$. The energy of the FD approximation to the solution 
was computed using centered second-order approximations to the derivatives 
for the inner points and one-sided second-order approximations to 
the spatial derivative for the vertices. The
integrals are approximated with the trapezoid rule. 
The time evolution of $E$ is shown in Fig. \ref{energy_mode1}. Though 
the energy for the exact solution is constant, the energy of the FD solution 
oscillates around a constant state, as the vertices produce 
small oscillations in time. These decay as $\Delta x^2$ with refinement.

\subsection{Discussion}

An important issue when solving a PDE on a metric graph 
is the choice of the number of modes (eigenvectors) for the representation
of the initial condition. Once the modes are fixed, the (linear)
PDE solution will evolve according to them. No new modes will be excited.

The number of modes necessary to resolve the initial condition
depends on its length scales. To illustrate this issue, consider
a Gaussian initial condition 
\be\label{gauss}
f(x) =   \exp(-{(x-x_0)^2 \over  s^2}), \ee
on arc $a$, of length $l_a$  with $x_0= l_a/2$ and  $s=l_a/ 10$. 
The other arcs are set to zero. We expand the initial condition
as in equation (\ref{expan2}).

\begin{figure}[H]
\centerline{
\epsfig{file=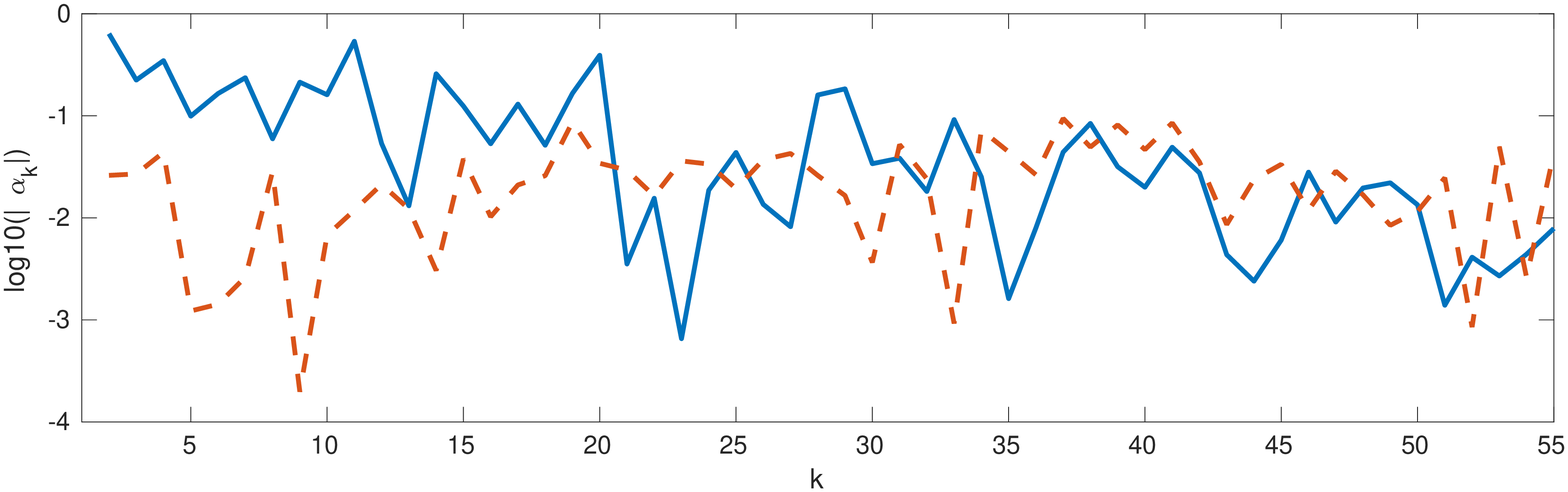,height=6cm,width=10 cm,angle=0}
}
\caption{Logarithms of the mode amplitudes $\alpha_q$ vs. $k_q$
for the
initial condition (\ref{gauss}) set on arc 1 ($l_1=11.91$) 
continuous line (blue) and arc 6 ($l_6=\sqrt{2}$) in dashed line (red). The 
parameters are $x_0= l/2$ and  $s=l/ 10$.}
\label{g14Gauss16}
\end{figure}

\begin{figure}[H]
\centerline{
\epsfig{file=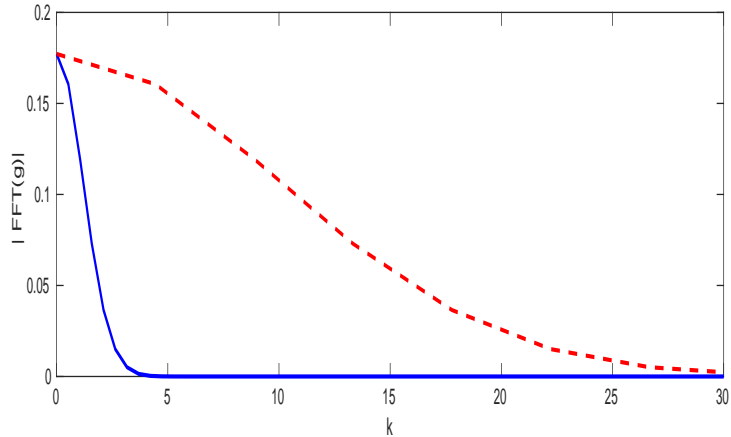,height=6cm,width=10 cm,angle=0}
}
\caption{Modulus of the one dimensional Fourier transform of the Gaussians
(\ref{gauss})
in arc 1 in continuous line (blue) and 6 in dashed line (red).}
\label{g14fourierGauss16}
\end{figure}

Fig. \ref{g14Gauss16} shows the logarithms of the mode amplitudes 
$\alpha_q$ vs. $k_q$
for the G14 graph for two Gaussian initial conditions
on arcs 1 ($l_1=11.91$) and 6 ($l_6=\sqrt{2}$). Clearly, the former 
is resolved as the amplitudes $|\alpha_q|$ decay. On the contrary,
for the initial condition set on arc 6, the $|\alpha_q|$s do not decay, so many more modes are needed.
This is confirmed by a one-dimensional Fourier analysis of the two Gaussians shown
in Fig. \ref{g14fourierGauss16}. The Gaussian on arc 1 extends up to
$k=5$ while the one on arc 6 extends up to $k=25$.
This simple analysis shows that a one-dimensional Fourier analysis provides
an estimate of the number of modes needed to resolve the initial condition.
For an inhomogeneous network such as G14, the number of
modes needed can be very large if the initial condition is concentrated
on a small arc.

\section{Conclusion}

In this article we developed and compared 
a spectral, a second-order finite difference and a discontinuous Galerkin 
method for solving linear PDEs on a metric graph with continuity and 
Kirchhoff vertex conditions.  
The spectral approach relies on a practical and robust 
algorithm for estimating the eigenvalues and eigenvectors
of the generalised Laplacian on a metric graph. It builds the matrix $M(k)$ and
finds the singular $k$s  using line optimisation on the rcond inverse condition number estimate and yields eigenvalues and eigenvectors of arbitrary order to machine precision.
We used this spectral formalism to solve a generalised Helmholtz problem, the Poisson equation, and the telegrapher's equation.

The FD method guarantees continuity of the solution as 
the computational node is placed at the vertices connecting the arcs. To 
achieve the second-order accuracy at the vertices, we placed a ghost 
grid point for each vertex. The standard second-order one-sided difference 
approximation to the fluxes in Kirchhoff's law reduced the accuracy of 
the solution at the node to the first-order as did the use of 
weighted arithmetic average of the nearby values.

The DG method allows us to include exact Kirchhoff and non-Kirchhoff 
flux conditions, but the continuity of the solution at 
the interval interfaces and the vertices requires the 
introduction of penalty terms.
This method provides arbitrary polynomial accuracy for spatial
approximations, but the condition number of the stiffness matrix 
grows exponentially, and for large problems it  requires the 
selection of an appropriate preconditioner.

The spectral method converges exponentially with the number of 
modes. It is superior for computing eigenvalues/eigenvectors
of the generalised Laplacian of arbitrary order, whereas the other 
two methods require the spatial resolution of highly oscillatory 
eigenvectors in such cases.
In addition, the spectral approach provides a very simple 
formalism for linear PDEs on graphs, exactly as for a one-dimensional linear PDE
on an interval. One first computes the spectral coefficients of the 
right hand side of the Poisson problem or the initial conditions 
for the telegrapher's equation.
Then the solution is written explicitly in terms of these modes.
Since the evolution problem is linear, each mode evolves 
separately, and no new modes appear.
{On the other hand, non-Kirchhoff vertex conditions, e.g. 
when the flux at each arc depends on the fluxes at the other arcs,  result 
in an overdetermined system and require special consideration to 
enforce continuity.
In this situation, the generalised Laplacian is not self-adjoint 
\cite{berkolaiko17} so the eigenvalues may become complex
and the eigenvectors not orthogonal. This
makes the spectral method much more complicated.}

Finally note that for some nonlinear problems, the spectral method
can be used as part of the time-split algorithms providing spectrally 
accurate solutions
to the linear part of the nonlinear equations, see \cite{spectralnls}.

\section*{Acknowledgements}
JGC acknowledges the support of the Agence Nationale de la recherche
through grant \textit{FRACTAL GRID}.
HK thanks the ARCS Foundation for support.

\appendix

\section*{Wave equation finite differences }
Here we describe a second-order finite difference method for approximating the solution to the wave equation on a network. We use the notation $u_{j,i}^r \approx u(x_j,t_r)$ on edge $e_i$.

For inner points $j=1, 2, \dots, N_i-1$ on edge $e_i$ we use a second-order approximation for the derivative in both time and space.

\be 0 = u_{xx} - u_{tt} = \frac{u_{j+1,i}^r - 2 u_{j,i}^r + u_{j-1,i}^r}{\Delta x_i^2} - \frac{u_{j,i}^{r+1} - 2 u_{j,i}^r + u_{j,i}^{r-1}}{\Delta t^2} + O(\Delta x_i^2 + \Delta t^2) \ee

The equation is solved for the current time-step $t=t_{r+1}$ and spatial step $x=x_{j}$ as follows
\be u_{j,i}^{r+1} = 2 u_{j,i}^r - u_{j,i}^{r-1} + \left( \frac{\Delta t}{\Delta x_i} \right)^2 \left[ u_{j+1,i}^r - 2 u_{j,i}^r + u_{j-1,i}^r \right] \ee

To implement the vertex conditions, we label arcs adjacent to the vertex as $c=1, 2, \dots, d$. To enforce the continuity condition we place a node exactly at the vertex and denote the solution value as $u_0$. The center vertex $u_0$, it is shared by each adjacent arc via the continuity condition $u_{c,0}=u_0 $. To implement the Kirchhoff flux condition, we use a centered, second-order scheme for the first derivative by extending each arc and adding a ghost point next to the vertex at $x=x_{-1,c}$ for each adjacent arc $e_c$. The derivative is taken in the outgoing direction from the center vertex.  
For $j=0$ we have 
\be 0 = \sum_{c=1}^d u_x(x_0) = \sum_{c=1}^d \frac{u_{1,c} - u_{-1,c}}{2 \Delta x_c}  + O(\Delta x_1^2 + \Delta x_2^2 + \cdots + \Delta x_d^2) \ee
Here we say $j=1$ is the point adjacent to the vertex regardless of edge orientation.

We then apply the finite difference scheme for the PDE. For the wave equation we have
\be \frac{u_{1,c}^r - 2 u_0^r + u_{-1,c}^r}{\Delta x_c^2} = \frac{u_0^{r+1} - 2 u_0^r + u_0^{r-1}}{\Delta t^2} \ee
This equation can be solved for $u_{-1,k}^r$ and substituted into the Kirchhoff flux equation to eliminate the ghost point.
\be 0 =  \sum_{c=1}^d \frac{u_{1,c}^r - \left[2 u_0^r - u_{1,c}^r + \left( \frac{\Delta x_c}{\Delta t}  \right)^2 \left[u_0^{r+1} - 2 u_0^r + u_0^{r-1} \right] \right]  }{2 \Delta x_c}   \ee

This equation can be solved for $u_0^{r+1}$ to find the solution at the vertex at time $t=t_{r+1}$.

\be u_0^{r+1} =  2 \Delta t^2 \frac{ \sum_{c=1}^d \frac{u_{1,c}^r}{\Delta x_c}  - u_0^r  \sum_{c=1}^d \frac{1}{\Delta x_c}}{ \sum_{c=1}^d \Delta x_c}+ 2 u_0^r - u_0^{r-1} \ee

\end{document}